\documentclass{svjour3}
\smartqed
\usepackage{latexsym,amsbsy,amssymb,amsmath}
\usepackage{color}
\usepackage{subcaption}
\usepackage{graphicx}
\usepackage{psfrag}
\usepackage{float}
\usepackage{t1enc}
\usepackage{amsfonts}
\usepackage{stmaryrd}
\usepackage{multirow}
\usepackage[table]{xcolor}

\usepackage{graphics}
\usepackage{epsf}
\usepackage{xfrac}
\usepackage{amscd}
\usepackage{wrapfig}
\usepackage[mathscr]{euscript}
\usepackage[english]{babel}
\usepackage{epsfig}
\usepackage{dsfont}
\everymath{\displaystyle}
\usepackage{setspace}
\usepackage{subcaption, algorithmicx}
\usepackage[]{algorithm}
\numberwithin{equation}{section}
\numberwithin{theorem}{section}

\newcommand{\bI}{\pmb{I}}
\newcommand{\bT}{\pmb{T}}

\newcommand{\bv}{\pmb{v}}
\newcommand{\bff}{\pmb{f}}

\newcommand{\bn}{\pmb{n}}
\newcommand{\bw}{\pmb{w}}

\newcommand{\bt}{\pmb{t}}

\newcommand{\bD}{\pmb{D}}

\newcommand{\bX}{\pmb{X}}

\newcommand{\bu}{\pmb{u}}

\newcommand{\mR}{\mathbb{R}}

\newcommand{\iS}{\mathcal{S}}
\newcommand{\iT}{\mathcal{T}}

\newcommand{\ahalf}{\sfrac{1}{2}}

\newcommand{\HL}[1]{\textcolor{black}{#1}}

\begin{document}

\title{A global-in-time domain decomposition method for the coupled nonlinear Stokes and Darcy flows 
\thanks{T.-T.-P. Hoang's work is partially supported by the US National Science Foundation under grant number DMS-1912626 and Auburn University's intramural grants program. \vspace{2pt}\\ 
H. Lee's work is partially supported by the US National Science Foundation under grant number DMS-1818842.}
}
\titlerunning{Global-in-time domain decomposition for the Stokes-Darcy coupling} % if too long for running head

\author{Thi-Thao-Phuong Hoang $\boldsymbol{\cdot} $ Hyesuk Lee
}

\authorrunning{T.-T.-P. Hoang and H. Lee} % if too long for running head

\institute{Thi-Thao-Phuong Hoang   \at
               Department of Mathematics and Statistics, Auburn University, Auburn, AL 36849. \\ 
              \email{tzh0059@auburn.edu}                  
           \and
             	Hyesuk Lee  \at
              School of Mathematical and Statistical Sciences, Clemson University, Clemson, SC 29634. \\
              \email{hklee@clemson.edu}
}

\date{Received: date / Accepted: date}

\maketitle

\begin{abstract}
We study a decoupling iterative algorithm based on domain decomposition for the time-dependent nonlinear Stokes-Darcy model, in which different time steps can be used in the flow region and in the porous medium. The coupled system is formulated as a space-time interface problem based on the interface condition for mass conservation. The nonlinear interface problem is then solved by a nested iteration approach which involves, at each Newton iteration, the solution of a linearized interface problem and, at each Krylov iteration, parallel solution of time-dependent linearized Stokes and Darcy problems. Consequently, local discretizations in time (and in space) can be used to efficiently handle multiphysics systems of coupled equations evolving at different temporal scales. Numerical results with nonconforming time grids are presented to illustrate the performance of the proposed method. 

\end{abstract}

%
% ------------------------------------------------
%
%    SECTION 1 : Introduction
%	
% ------------------------------------------------
%
\section{Introduction}
\label{Sec:intro}  
Multiscale and multiphysics processes are ubiquitous in many science and engineering applications. Mathematically, coupled partial differential equations are used to model various processes possibly taking place on different regions of the problem domain and at different scales in space and time. 
%The spatial and temporal scales associated with various processes may be vary with several orders of magnitude. 
One example of such a coupling is the coupled (Navier-)Stokes-Darcy system arising in a number of applications: surface and subsurface flow interaction, flow in vuggy porous media, industrial filtrations, biofluid-organ interaction, cardiovascular flows, and others. In these applications, the Stokes equations are used to model the free flow and the Darcy equations are used to model the flow in a porous medium; the two flow domains are coupled via suitable transmission conditions on the interface to enforce mass conservation, balance of the normal forces and the Beavers-Joseph-Saffman law \cite{BJ67,S71,JM96}.

The development of numerical approximations and efficient solvers for the Stokes-Darcy coupling has been an active research area and attracted great attention over the past two decades. For the stationary case, the existence and uniqueness of the weak solution of the coupled system are proved in \cite{Disca02,Layton03,Bernardi08,Disca09}. %for steady state and time-dependent cases. 
Regarding a numerical solution of the mixed Stokes-Darcy model, one can either solve the coupled system directly with some suitable preconditioner, or use the domain decomposition-based approach \cite{QV99,Toselli:DDM:2005} to decouple the system into two local subsystems which are solved separately. Concerning the former or monolithic approach, new finite element spaces were studied in \cite{Mardal02,AB07,Burman07,AG09} with mixed formulations and in \cite{Riviere05,RY05} with discontinuous approximations. Preconditioning techniques for solving the sparse linear system of saddle point form resulted from finite element discretization of the fully coupled Stokes-Darcy system were investigated in \cite{Mu09,Marquez13,CLS16}. Concerning the decoupled approach, several directions have been considered. Lagrange multiplier techniques were proposed in \cite{Layton03,Ervin09} and mortar finite elements were studied in \cite{GS07,Bernardi08,Ervin11,Girault14,Yotov17} in which the meshes on the interface and subregions do not necessarily match. Heterogeneous domain decomposition methods were explored using either the classical Dirichlet-Neumann (Steklov-Poincar\'e) type operator  \cite{Disca02,Disca03,Disca04,Hoppe07,GS10,VWY14} or the Robin-Robin interface conditions \cite{Disca07,Disca18,Cao11,Chen11,Caia14}. Two-grid methods were applied to the mixed Stokes-Darcy model in \cite{Mu07,Mu12}, and optimization-based approach was proposed in \cite{Ervin14}. 

For nonstationary Stokes-Darcy problems, only a few studies have been carried out. A monolithic method based on implicit time discretization was presented in~\cite{DiscaThesis} in which the evolutionary system is uncoupled at each time step by domain decomposition iteration. In \cite{Mu10}, a decoupled backward Euler scheme was devised by lagging the interface coupling terms, i.e. at each time level, one solves the Stokes and Darcy problems using Neuman interface boundary conditions computed from the previous time level. Long term stability of this method and a modified two-step method was analyzed in \cite{Layton13}. A similar decoupled scheme with Robin interface conditions was studied in \cite{Cao14} in which higher order time discretization (three-step backward differentiation method) was also considered. In these works, the same time step is used in both regions. Decoupled schemes with different time step sizes were proposed and analyzed in~\cite{Shan13,Rybak14}. These schemes are extensions of the method in \cite{Mu10} in which the time step size in the Stokes region is an integral multiple of the time step size in the Darcy region. The advancement in time is then carried out sequentially; first the Stokes problem is solved with a small time step size using the Darcy pressure (freezing from the previous coarse time step) as interface data, then the Darcy problem is solved using the recently computed Stokes velocity as interface data. These methods are non-iterative by using an explicit method for the coupling terms, and the key issue is how to achieve desired accuracy and stability properties. A different approach was proposed  in \cite{HSLee14} by formulating the coupled problem as a constrained optimal control problem which is solved at each time step by a least square method (thus, the same time step is used in both regions). 

As the model concerns the flow of fluid, there are two possible fluid types: Newtonian fluids (e.g. water and air) and non-Newtonian fluids (e.g. honey and quicksand). The difference between these two types of fluid lies in the viscosity which is a constant for Newtonian fluids and a function of the magnitude of the deformation tensor for non-Newtonian fluids (more discussion can be found in \cite{Ervin09}). Mathematically, one deals with a linear or nonlinear coupled flow problem; the nonlinear Stokes-Darcy coupling was considered in \cite{Ervin09,Ervin11,HSLee14,Ervin14}. In addition, approximation methods for the nonlinear Navier-Stokes/Darcy system were studied in \cite{Disca07,Girault09,Mu09,CR09,CR10,CesR09}, and for the coupling with transport in \cite{VY09,CGR13,RZ17}. 

In this work, we aim to develop a parallel decoupling method for the time-dependent nonlinear Stokes-Darcy system in which different time step sizes can be used in the free flow domain and the porous medium. Differently from~ \cite{Shan13,Rybak14}, we apply the so-called global-in-time (or space-time) domain decomposition method in which the dynamic system is decoupled into dynamic subsystems defined on the subdomains (resulting from a spatial decomposition), then time-dependent problems are solved in each subdomain at each iteration and information is exchanged over space-time interfaces between subdomains. Consequently, local discretizations in both space and time can be enforced in different regions of the computational domain, which makes the method well-suited and efficient for multiscale multiphysics problems. Note that this approach is implicit in time, thus considerably large time step sizes can be used without affecting stability, unlike the explicit method in \cite{Shan13,Rybak14}. This can be important for applications in geosciences where long time simulations are often required. The method has been studied for porous medium flows (see \cite{H13,H16} and the references therein), and here we extend the idea to the nonlinear Stokes-Darcy coupling, which, to the best of our knowledge, hasn't been considered in the literature. We construct a time-dependent Steklov-Poincar\'e type operator and reduce the coupled problem into a nonlinear time-dependent interface problem. The interface problem is then solved by a nested iteration approach which involves, at each Newton iteration, the solution of a linearized interface problem and, at each Krylov iteration, parallel solution of time-dependent linearized Stokes and Darcy problems. As the local problems are solved globally in time at each iteration, it makes possible the use of different time discretization methods or different time grids in the Stokes and Darcy regions.  To exchange information at the interface with nonconforming time grids, an $L^{2}$ time projection between subdomains is performed by an optimal projection algorithm without any additional grid \cite{Gander05}. High order time stepping methods can be applied straightforwardly, see \cite{Japhet12}. The idea can be generalized to the case of multiple subdomains where interfaces of different types are introduced: Stokes-Darcy, Stokes-Stokes and Darcy-Darcy as considered in \cite{VWY14} for the steady problems. However, in this work we restrict ourselves to the case of two subdomains and conforming spatial meshes, and focus on numerical performance -  in terms of accuracy and efficiency - of the proposed method with nonmatching time grids.

The rest of this paper is structured as follows. In Section~\ref{sec:model}, we present the model problem which is the nonstationary nonlinear Stokes-Darcy system, and the interface coupling conditions. The variational formulation of the continuous coupled system is derived in Section~\ref{sec:weakform}. The coupled problem is formulated as a time-dependent nonlinear interface problem in Section~\ref{sec:IP}, and nonconforming time discretization is discussed in Section~\ref{sec:time}. Numerical results are presented in Section~\ref{sec:NumRe} to study the performance of the proposed algorithm with nonmatching time grids. Finally, some concluding remarks are given in Section~\ref{sec:conclu}.
%
% ------------------------------------------------
%
%    SECTION: Model problem
%	
% ------------------------------------------------
%
\section{Time-dependent nonlinear Stokes-Darcy system}
\label{sec:model}
\begin{figure}[ht!]
\centering
\includegraphics[scale=0.4]{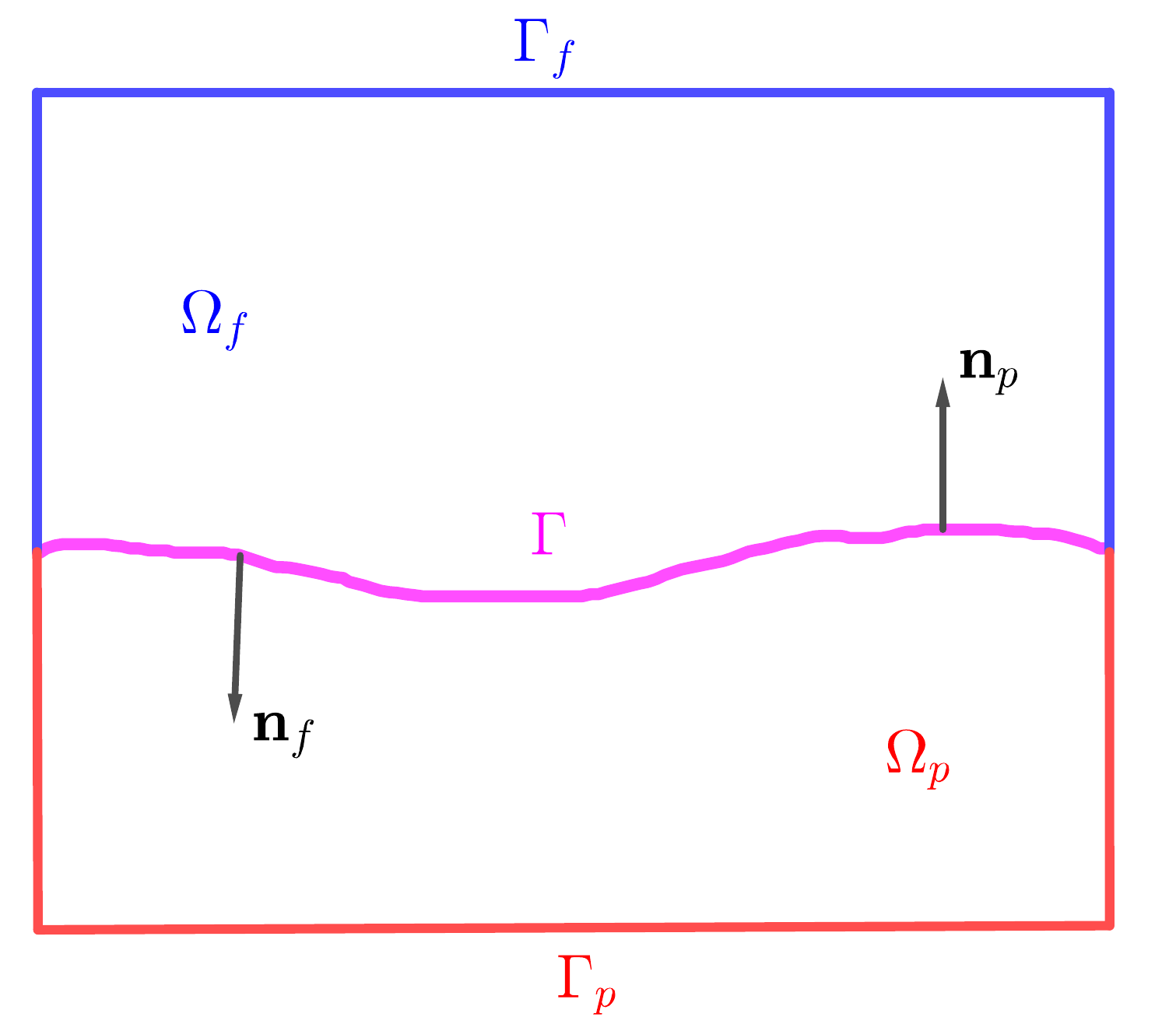}
\caption{Example of a two dimensional domain formed by a fluid region and a porous medium.}
\label{fig:domain}
\end{figure}

We consider a free non-Newtonian fluid flow in $\Omega_{f}$ coupled with a porous medium flow in $\Omega_{p}$, where $\Omega_{f}$ and $\Omega_{p}$ are subsets of $\mR^{d}$ for $d=2,3$. Denote by $\Gamma$ the interface between the two domains, and by $\Gamma_{f}= \partial \Omega_{f} \setminus \Gamma$ and $\Gamma_{p}=\partial \Omega_{p} \setminus \Gamma$ the external boundaries of the fluid domain and porous medium respectively (see Figure~\ref{fig:domain}). Let $\bn_{f}$ and $\bn_{p}$ be the outward unit normal vectors to $\Omega_{f}$ and $\Omega_{p}$ respectively, and $\{ \bt_{j} \}_{j=1, \hdots, d-1}$ be an orthogonal set of unit tangent vectors on $\Gamma$.\vspace{3pt}

Let $T >0$ be a finite time. The free flow in $\Omega_{f}$ is described by the nonlinear Stokes equations subject to no-slip boundary condition on $\Gamma_{f}$:
\begin{subequations}
\label{eq:Stokes}
\begin{align}
\frac{\partial \bu_{f}}{\partial t}-\nabla\cdot\bT(\bu_{f}, p_{f})
%=\frac{\partial \bu_{f}}{\partial t} -\nabla \cdot \left (\nu_{f}(\vert \bD(\bu_{f}) \vert) \bD(\bu_{f})\right ) + \nabla p_{f} 
& = \bff_{f} & \text{in} \; \Omega_{f} \times (0,T), \\
\nabla\cdot \bu_{f} & =0 & \text{in} \; \Omega_{f} \times (0,T), \\
\bu_{f} & = \pmb{0} & \text{on} \; \Gamma_{f} \times (0,T), \\
\bu_{f}(\cdot, 0) & = \bu_{f0} & \text{in} \; \Omega_{f},
\end{align}
\end{subequations}
where $\bu_{f}$ is the fluid velocity, $p_{f}$ the fluid pressure, $\bT(\bu_{f}, p_{f})=\nu_{f}(\vert \bD(\bu_{f}) \vert) \bD(\bu_{f})-p_{f}\bI$ the stress tensor (with $\bI$ the identity tensor), $\bD(\bu_{f}) = \frac{1}{2} \left (\nabla \bu_{f} + \nabla \bu_{f}^{T}\right )$ the rate of the strain tensor, $\nu_{f}(\cdot)$ the fluid viscosity and $\bff_{f}$  the body force. In this work we consider the Cross model for the viscosity function:
\begin{equation} \label{eq:visf}
\nu_{f}(\vert \bD(\bu_{f}) \vert)=\nu_{f\infty} +\frac{\nu_{f0}-\nu_{f\infty}}{1+K_{f} \vert \bD(\bu_{f})\vert^{2-r_f}},
\end{equation}
where $r_{f}>1$, $\nu_{f\infty}, \nu_{f0} >0$ and $K_{f}>0$ are constants; $\nu_{f\infty}$ and $\nu_{f0}$ denote the limiting viscosity values at an infinite shear rate and at zero shear rate respectively, and satisfy $\nu_{f \infty} \leq \nu_{f0}$. Other nonlinear viscosity models such as Carreau model, power law model and Ladyzhenskaya model can also be used \cite{Ervin09}. 

The porous medium flow in $\Omega_{p}$ is described by the nonlinear Darcy equations subject to no-flux boundary condition on $\Gamma_{p}$:
\begin{subequations}
\label{eq:Darcy}
\begin{align}
{\nu_{\text{eff}}(\vert \bu_{p}\vert)}  \, \kappa^{-1} \bu_{p} + \nabla p_{p} & = 0 &\text{in} \; \Omega_{p} \times (0,T), \\
S_{p} \frac{\partial p_{p}}{\partial t} + \nabla \cdot \bu_{p} & = f_{p} & \text{in} \; \Omega_{p} \times (0,T), \\
\bu_{p} \cdot \bn_{p} & = 0 &\text{on} \; \Gamma_{p} \times (0,T), \\
p_{p}(\cdot, 0) & = p_{p0} & \text{in} \; \Omega_{p},
\end{align}
\end{subequations}
where  $\bu_{p}$ and $p_{p}$ are the Darcy velocity and pressure respectively, $S_{p} > 0$ the storage coefficient, ${\nu_{\text{eff}}(\cdot)}$ the effective fluid viscosity, $\kappa >0 $ the permeability and $f_{p}$ is the source/sink. The Cross model for $\nu_{\text{eff}}$ is defined as follows (see \cite{Ervin09} for other models):
\begin{equation} \label{eq:visp}
\nu_{p}(\vert \bu_{p}\vert)= \nu_{p\infty} +\frac{\nu_{p0}-\nu_{p\infty}}{1+K_{p} \vert \bu_{p}\vert^{2-r_p}},
\end{equation}
where $r_{p}>1$, $ \nu_{p0} \geq \nu_{p\infty} > 0$ and $K_{p}>0$ are constants. 

The viscosity functions \eqref{eq:visf} and \eqref{eq:visp} have the following properties which will be used in later analysis \cite{Ervin09}: \vspace{-0.2cm}
\begin{itemize}
\item[(A1)] $\nu_{f}(\cdot)$ and $\nu_{p}(\cdot)$ are strongly monotone and bounded from below and above by positive constants. 
%$$  \nu_{f\infty} \leq \nu_{f}(\vert \bu \vert) \leq \nu_{f0}, \quad \nu_{p\infty} \leq \nu_{p}(\vert \bu \vert) \leq \nu_{p0}, \quad \forall \bu \in \mR^{d}.$$
\item[(A2)] The nonlinear functions $\nu_{f}(\vert \bu \vert) \vert \bu \vert$ and $\nu_{p}(\vert \bu \vert) \bu$ are uniformly continuous with respect to $\bu \in \mR^{d}$.  \vspace{-0.2cm}
\end{itemize}
Note that the standard linear Stokes-Darcy system can be recovered by setting $r_{f}=r_{p}=2$. 

The coupled Stokes-Darcy system is closed by the following coupling conditions on the space-time interface:
\begin{subequations}
\label{eq:TCs}
\begin{align}
\bu_{f} \cdot \bn_{f} + \bu_{p} \cdot \bn_{p} & = 0 \quad \text{on} \; \Gamma \times (0,T), \label{eq:TC1}\\
-\bn_{f} \cdot (\nu_{f}(\vert \bD(\bu_{f})\vert) \bD(\bu_{f})-p_{f}\bI) \cdot \bn_{f} & = p_{p} \quad \text{on} \; \Gamma \times (0,T), \label{eq:TC2}\\
-\bn_{f} \cdot (\nu_{f}(\vert \bD(\bu_{f})\vert) \bD(\bu_{f})-p_{f}\bI ) \cdot \bt_{j} & = c_{BJS} \bu_{f} \cdot \bt_{j}\quad  \text{on} \; \Gamma \times (0,T), \; j=1, \hdots, d-1,\label{eq:BJS}
 \end{align}
 \end{subequations}
where $c_{BJS}$ is a positive constant. These coupling conditions have been studied extensively in the literature (e.g. \cite{Disca02,Layton03,Cao10a}). The first two conditions enforce the continuity of the normal component of velocities and the continuity of the normal stress respectively. The third condition is the Beavers-Joseph-Saffmann condition \cite{S71,JM96}, stating the connection between the slip velocity and the shear stress along the interface. It is a simplification of the Beavers-Joseph condition \cite{BJ67} by neglecting the porous medium velocity tangent to the interface. Thus \eqref{eq:BJS} is actually not a coupling condition as it only involves the fluid domain's variables. Next, we derive the weak formulation of the coupled system with the use of Lagrange multipliers.

%
% ------------------------------------------------
%
%    SECTION: Weak form of fully coupled system
%	
% ------------------------------------------------
%
\section{Variational formulation of the fully coupled system} \label{sec:weakform}
In the following, we will use the convention that if $V$ is a space of
functions, then we write $\pmb{V}$ for a space of vector functions having each component
in $V$. In order to write the variational formulation of the coupled problems, we first introduce the functional spaces:
\begin{align*}
&\bu_{f} \in \bX_{f}:=\{ \bv \in \pmb{H^{1}}(\Omega_{f}): \bv = \pmb{0} \; \text{on} \; \Gamma_{f} \},  \quad p_{f} \in Q_{f}:=L^{2}(\Omega_{f}), \\
&\bu_{p} \in \bX_{p}:=\{ \bv \in \pmb{L^{2}}(\Omega_{p}): \nabla \cdot \bv \in L^{2}(\Omega_{p}), \; \bv \cdot \bn_{p} = 0 \; \text{on} \; \Gamma_{p}\}, \quad p_{p} \in Q_{p}:=L^{2}(\Omega_{p}).
\end{align*}
Let $\Omega=\Omega_{1}  \cup \Omega_{2} \cup \Gamma$ and define the spaces $\bX$ and $Q$ on $\Omega$ by $ \bX = \bX_{f} \times \bX_{p}$ and $Q=Q_{f} \times Q_{p}$ respectively. 
Denote by $\bX_{f}^{\ast}$ the dual space of $\bX_{f}$. 
For a domain $\Theta = \Omega_{f}$ or $\Theta = \Omega_{p}$, we denote by $(\cdot, \cdot)_{\Theta}$ the $L^{2}$ inner product over $\Theta$. % when the context is clear, the subscript can be omitted. 
%Let $H^{\ahalf}(\Gamma)$ be the trace space of $H^{1}(\Omega_{f})$ or $H^{1}(\Omega_{p})$ on $\Gamma$, and $H^{-\ahalf}(\Gamma)$ be the dual space of $H^{\ahalf}(\Gamma)$. Denote by $\wH^{\ahalf}(\Gamma)$ the closure of the trace of $\bX_{f} \cdot \bn_{f}$ on $\Gamma$ with respect to the $H^{\ahalf}(\Gamma)$ norm, then for $\Gamma$ sufficiently smooth, $\wH^{\ahalf}(\Gamma) \subset H^{\ahalf}(\Gamma)$. The dual space of $\wH^{\ahalf}(\Gamma)$ is $\wH^{-\ahalf}(\Gamma)$ \cite{Ervin11}. 
As in the stationary case~\cite{Layton03,Ervin09}, we introduce the Lagrange multiplier $\lambda$ on the interface representing:
\begin{equation} \label{eq:lambda}
\lambda = -\bn_{f} \cdot (\nu_{f}(\vert \bD(\bu_{f})\vert) \bD(\bu_{f})-p_{f}\bI) \cdot \bn_{f}  = p_{p} \quad \text{on} \; \Gamma \times (0,T).
\end{equation}
\HL{The space for the Lagrange multiplier is $\Lambda:=H^{\ahalf}_{00}(\Gamma)$ (see~\cite{Layton03}). We denote by $\Lambda^{\ast}:=\left (H^{\ahalf}_{00}(\Gamma)\right )^{\ast}$ the dual space of $\Lambda$ and by $\langle \cdot, \cdot \rangle_{\Gamma}$ the duality pairing between $\Lambda^{\ast}$ and $\Lambda$.}
Define the bilinear forms $a(\cdot, \cdot): \bX \times \bX \rightarrow \mR$, $b(\cdot, \cdot): \bX \times Q  \rightarrow \mR$ and $b_{I}: \bX \times \Lambda \rightarrow \mR$ by:
\begin{align*}
a(\bu,\bv) &= a_{f}(\bu_{f}, \bv_{f})+ a_{p}(\bu_{p}, \bv_{p}), \quad  b(\bv, q)  =b_{f}(\bv_{f},q_{f}) + b_{p}(\bv_{p}, q_{p}),\\  
b_{\Gamma}(\bv, \zeta)  &= b_{\Gamma f}(\bv_{f}, \zeta) + b_{\Gamma p}(\bv_{p}, \zeta), 
\end{align*}
where
\begin{align*}
a_{f}(\bu_{f}, \bv_{f})&=  \left (\nu_{f}(\vert \bD(\bu_{f})\vert) \bD(\bu_{f}), \bD(\bv_{f})\right )_{\Omega_{f}}  + \sum_{j=1}^{d-1} c_{BJS}(\bu_{f} \cdot \bt_{j}, \bv_{f} \cdot \bt_{j})_{\Gamma}, \\
a_{p}(\bu_{p}, \bv_{p})&=\left ({\nu_{\text{eff}}(\vert \bu_{p}\vert)}  \, \kappa^{-1}  \bu_{p}, \bv_{p}\right )_{\Omega_{p}} ,\\
b_{f}(\bv_{f}, q_{f}) &=(q_{f}, \nabla \cdot \bv_{f})_{\Omega_{f}}, \quad
b_{p}(\bv_{p}, q_{p}) = (q_{p}, \nabla \cdot \bv_{p})_{\Omega_{p}}, \\
 b_{ \Gamma f}(\bv_{f}, \zeta) &=\langle \zeta, \bv_{f} \cdot \bn_{f} \rangle_{\Gamma}, \quad b_{\Gamma p}(\bv_{p}, \zeta)= \langle \bv_{p} \cdot \bn_{p}, \zeta \rangle_{\Gamma}.
\end{align*}
%Note that $\langle \zeta, \bv_{f} \cdot \bn_{f} \rangle_{\Gamma}$ denotes a duality pairing between $\wH^{-\ahalf}(\Gamma)$ and $\wH^{\ahalf}(\Gamma)$, while $\langle \bv_{p} \cdot \bn_{p}, \zeta \rangle_{\Gamma}$ denotes a duality pairing between $H^{-\ahalf}_{00}(\Gamma)$ and $H^{\ahalf}_{00}(\Gamma)$.
%
The weak formulation of the coupled system \eqref{eq:Stokes}-\eqref{eq:Darcy}-\eqref{eq:TCs} is then written as follows (detailed derivation for the stationary problems can be found in \cite{Ervin09}): \vspace{0.2cm}\\
%{\em Given $\bff_{f} \in L^{2}(0,T;\bX_{f}^{\ast})$, $f_{p} \in L^{2}(0,T; Q_{p})$, $\bu_{f0} \in \pmb{L^{2}}(\Omega_{f})$ and $p_{p0} \in Q_{p}$,
{\em For a.e. $t \in (0,T)$, find $\left (\bu(t), p(t),\lambda(t)\right ) \in \bX \times Q \times \HL{\Lambda}$ such that:}
\begin{align}
(\partial_{t} \bu_{f}, \bv_{f})_{\Omega_{f}} + a(\bu, \bv) - b(\bv,p) + b_{\Gamma}(\bv, \lambda) &= (\bff_{f}, \bv_{f})_{\Omega_{f}}, & \forall \bv \in \bX, \label{eq:weak1}\\
 b(\bu,q) - b_{\Gamma}(\bu, \zeta)+(S_{p}\partial_{t} p_{p}, q_{p})_{\Omega_{p}} &= (f_{p},q_{p})_{\Omega_{p}}, & \forall (q,\zeta) \in Q \times \HL{\Lambda}, \label{eq:weak2} 
\end{align}
{ \em with the initial conditions}
$$ \bu_{f}(\cdot, 0) = \bu_{f0} \quad \text{in} \; \Omega_{f}, \qquad p_{p} (\cdot, 0) = p_{p0} \quad \text{in} \; \Omega_{p}.
$$

The existence and uniqueness of the weak solution to the non-stationary and linear Stokes-Darcy system is proved in \cite{Cao10a} using the Stokes-Laplace formulation, i.e. the velocity and pressure are the unknowns in the fluid flow domain and the pressure is the only unknown in the porous media domain. In addition, no Lagrange multiplier is introduced and the physically more accurate coupling condition - the Beavers-Joseph condition - is considered in \cite{Cao10a}. The well-posedness of the stationary nonlinear Stokes-Darcy system in mixed form with a Lagrange multiplier is proved in \cite{Ervin09}. Here we assume the variational formulation \eqref{eq:weak1}-\eqref{eq:weak2} is well-posed, and focus on the decoupled approach based on global-in-time domain decomposition. 

\section{Decoupled problems and nested iteration approach}  \label{sec:IP}
We shall reformulate the Stokes-Darcy coupled problem as a space-time interface problem with the interface unknown $\lambda$ defined in \eqref{eq:lambda}. Assume that $\lambda$ is given, the Stokes and Darcy problems are then decoupled. We derive the weak formulations of the local problems using \eqref{eq:lambda} as boundary conditions on the interface, then formulate the interface problem which is solved by a nested iteration approach. 
\subsection{Free fluid flow}
We first consider the Stokes problem with Neumann boundary condition on the interface~$\Gamma$:
\begin{equation} \label{eq:lamS}
 -\bn_{f} \cdot (\nu_{f}(\vert \bD(\bu_{f})\vert) \bD(\bu_{f})) \cdot \bn_{f} + p_{f}=\lambda, \quad \text{on} \; \Gamma \times (0,T).
 \end{equation}
Its variational formulation is given by: \\
{\em For a.e. $t \in (0,T)$, find $\left (\bu_{f}(t), p_{f}(t)\right ) \in \bX_{f} \times Q_{f}$ such that:}
\begin{align}
\left (\partial_{t} \bu_{f}, \bv_{f}\right ) + a_{f}(\bu_{f}, \bv_{f}) - b_{f}(\bv_{f},p_{f}) &= (\bff_{f}, \bv_{f})_{\Omega_{f}}- b_{\Gamma f}(\bv_{f}, \lambda), & \forall \bv_{f} \in \bX_{f}, \vspace{0.1cm} \label{eq:Stokes1w}\\
b_{f}(\bu_{f},q_{f}) &= 0, & \forall q_{f} \in Q_{f}, \label{eq:Stokes2w} 
\end{align}
{\em with the initial condition}
\begin{equation}
\bu_{f}(\cdot,0) = \bu_{f0}, \quad \text{in} \; \Omega_{f}.   \label{eq:Stokesic}
\end{equation}

For given $\bff_{f} \in L^{2}(0,T; \bX_{f}^{\ast})$, $\lambda \in L^{2}(0,T; \HL{\Lambda})$ and $\bu_{f0} \in \bX_{f}$, the existence and uniqueness of the solution 
$$(\bu_{f}, p_{f}) \in \left (H^{1}(0,T; \pmb{L^{2}}(\Omega_{f}) \cap L^{2}(0,T; \bX_{f})\right ) \times L^{2}(0,T; Q_{f})$$
 to \eqref{eq:Stokes1w}-\eqref{eq:Stokes2w} with the initial condition \eqref{eq:Stokesic} are followed from the strong monotonicity of the viscosity function \eqref{eq:visf}, \cite{Ervin11} and the classical result of wellposedness of evolutionary (Navier-)Stokes equations \cite[Chapter III]{Temam77}.

\subsection{Porous medium flow}
We now consider the Darcy flow with Dirichlet boundary condition on the interface $\Gamma$:
\begin{equation} \label{eq:lamD}
p_{p} = \lambda, \quad \text{on} \; \Gamma \times (0,T).
\end{equation}
Its variational formulation is given by: \\
{\em For a.e. $t \in (0,T)$, find $\left (\bu_{p}(t), p_{p}(t)\right ) \in \bX_{p} \times Q_{p}$ such that:}
\begin{align}
a_{p}(\bu_{p}, \bv_{p}) - b_{p}(\bv_{p},p_{p}) &= - b_{\Gamma p}(\bv_{p}, \lambda), &\forall \bv_{p} \in \bX_{p}, \label{eq:Darcy1w}\\
b_{p}(\bu_{p}, q_{p})+\left (S_{p} \partial_{t} p_{p}, q_{p}\right ) &= (f_{p}, q_{p})_{\Omega_{p}}, & \forall q_{p} \in Q_{p}, \label{eq:Darcy2w}
\end{align}
{\em with the initial condition}
\begin{equation}
p_{p}(\cdot,0) = p_{p0}, \quad \text{in} \; \Omega_{p}.   \label{eq:Darcyic}
\end{equation}

For given $f_{p} \in L^{2}(0,T; Q_{p}))$, $\lambda \in L^{2}(0,T; \HL{\Lambda})$ and $p_{p0} \in H^{1}(\Omega_{p})$, there exists a unique solution 
$$(\bu_{p}, p_{p}) \in L^{2}(0,T; \bX_{p}) \times H^{1}(0,T; Q_{p})$$
 to \eqref{eq:Darcy1w}-\eqref{eq:Darcy2w} with the initial condition \eqref{eq:Darcyic}. This is obtained by using the strong monotonicity of the viscosity function \eqref{eq:visp} and the Faedo-Galerkin method for mixed formulations of the Darcy problem as in \cite{H13}.

\subsection{Nonlinear space-time interface problem}
We first introduce the interface operators:
$$
\begin{array}{lcccl}
	\iS_f\,:& L^{2}(0,T; \HL{\Lambda}) & \longrightarrow & L^{2}(0,T;\HL{\Lambda^{\ast}}), & S_{f}(\lambda)=\bu_{f}(\lambda) \cdot \bn_{f} \vert_{\Gamma}, \vspace{0.2cm}\\
	\iS_p\,:& L^{2}(0,T; \HL{\Lambda})  & \longrightarrow & L^{2}(0,T;\HL{\Lambda^{\ast}}), & S_{p}(\lambda)=\bu_{p}(\lambda) \cdot \bn_{p} \vert_{\Gamma},
\end{array}
$$
where $\left (\bu_{f}(\lambda),p_{f}(\lambda)\right )$ and  $\left (\bu_{p}(\lambda),p_{p}(\lambda)\right )$ are the solutions to the Stokes problem \eqref{eq:Stokes1w}-\eqref{eq:Stokesic} and the Darcy problem \eqref{eq:Darcy1w}-\eqref{eq:Darcyic} respectively. 

As the continuity of the normal stress \eqref{eq:TC2} is imposed via $\lambda$ in \eqref{eq:lamS} and \eqref{eq:lamD} (note that the Beaver-Joseph-Saffmann condition is imposed naturally in \eqref{eq:Stokes1w}), there remains to enforce the condition \eqref{eq:TC1}, which leads to the interface problem: \\
{\em For a.e. $t \in (0,T)$, find $\lambda(t) \in \HL{\Lambda}$ such that:}
\begin{equation} \label{eq:nonIF}
\int_{0}^{T} \biggl (\langle S_{f}(\lambda), \zeta \rangle_{\Gamma} + \langle S_{p} (\lambda), \zeta\rangle_{\Gamma} \biggl ) \, ds= 0, \quad \forall \zeta \in \HL{\Lambda}.
\end{equation}
This is a time-dependent and nonlinear problem which will be solved by a nested iteration approach. Toward that end, we define the operator:
\begin{equation}
\Psi (\lambda) := S_{f}(\lambda)+ S_{p} (\lambda),
\end{equation}
and apply the Newton algorithm to \eqref{eq:nonIF} to obtain the following linear system at each iteration $k$:
\begin{equation} \label{eq:linIF}
\int_{0}^{T} \left \langle \Psi^{\prime}(\lambda^{k})(\lambda^{k+1}-\lambda^{k}), \zeta \right \rangle_{\Gamma} \, ds = \int_{0}^{T}\left  \langle -\Psi(\lambda^{k}),  \zeta \right  \rangle_{\Gamma} \, ds,  \quad \forall \zeta \in \Lambda,
\end{equation}
where  $\Psi^{\prime}(\lambda)(h) = S_{f,\lambda}^{\text{lin}}(h)+ S_{p, \lambda}^{\text{lin}}(h)$, and
\begin{align*}
& S_{f,\lambda}^{\text{lin}}(h) = \bw_{f}(h) \cdot \bn_{f} \vert_{\Gamma}, \quad  S_{p, \lambda}^{\text{lin}}(h) = \bw_{p}(h)\cdot \bn_{p})\vert_{\Gamma},
\end{align*} 
in which $\left (\bw_{f}(h), \xi_{f}(h)\right )$ is the solution to the linearized Stokes problem \cite{HSLee14}:
\begin{align}
&\left (\partial_{t} \bw_{f}, \bv_{f}\right ) + \left (\nu_{f}(\vert \bD(\bu_{f})\vert) \bD(\bw_{f}), \bD(\bv_{f})\right ) \nonumber \\
& + \left ( \frac{(r_{f}-2)(\nu_{f0}-\nu_{f\infty})K_{f}}{(1+K_{f} \vert \bD(\bu_{f})\vert^{2-r_{f}})^{2} \vert \bD(\bu_{f})\vert^{r_{f}}} \bD(\bu_{f}) (\bD(\bu_{f}):\bD(\bw_{f})), \bD(\bv_{f}) \right ) \nonumber \\
&- (\xi_{f}, \nabla \cdot \bv_{f}) + \sum_{j=1}^{d-1} c_{BJS}(\bw_{f} \cdot \bt_{j}, \bv_{f} \cdot \bt_{j})_{\Gamma}= -\langle h, \bv_{f} \cdot \bn_{f} \rangle_{\Gamma}, \quad \forall \bv_{f} \in \bX_{f}, \vspace{0.1cm} \label{eq:linStokes1}\\
& (q_{f}, \nabla \cdot \bw_{f}) = 0, \quad \forall q_{f} \in Q_{f}, 
\end{align}
and  $\left (\bw_{p}(h), \xi_{p}(h)\right )$ is the solution to the linearized Darcy problem:
\begin{align}
& \left (S_{p} \partial_{t} \xi_{p}, q_{p}\right ) + (q_{p}, \nabla \cdot \bw_{p}) = 0, \quad \forall q_{p} \in Q_{p}, \\
& \left ({\nu_{\text{eff}}(\vert \bu_{p}\vert)}  \, \kappa^{-1}  \bw_{p}, \bv_{p}\right ) + \left ( \frac{(r_{p}-2)(\nu_{p0}-\nu_{p\infty})K_{p}}{(1+K_{p} \vert\bu_{p}\vert^{2-r_{p}})^{2} \vert \bu_{p} \vert^{r_{p}}} \bu_{p} (\bu_{p}:\bw_{p}), \bv_{p} \right ) - (\xi_{p}, \nabla \cdot \bv_{p}) \nonumber \\
& \hspace{4cm} = -\langle h, \bv_{p} \cdot \bn_{p} \rangle_{\Gamma}, \quad \forall \bv_{p} \in \bX_{p},\label{eq:linDarcy2}
\end{align}
Note that $\bu_{f}=\bu_{f}(\lambda)$ in \eqref{eq:linStokes1} and $\bu_{p}=\bu_{p}(\lambda)$ in \eqref{eq:linDarcy2}. The nested iteration algorithm for solving \eqref{eq:nonIF} is summarized in Algorithm~1. \vspace{0.2cm}
		\begin{flushleft}
		{\bf Algorithm 1 - Nested Iteration Approach} \vspace{0.1cm}\\
		\textbf{Input}: $\lambda^{0}$ initial guess, $\epsilon$ tolerance and $N_{\text{iter}}$ maximum number of iterations. \\
		\textbf{Output}: $\lambda^{k}$ \vspace{0.2cm}\\
		$k=0$, $\text{error} =0,$\\
		\textbf{while} $k < N_{\text{iter}}$ and $\text{error}> \epsilon$, \textbf{do}:
		\end{flushleft}
		\begin{algorithmic}[1]
			\item Compute the RHS of \eqref{eq:linIF} by solving the nonlinear Stokes problem \eqref{eq:Stokes1w}-\eqref{eq:Stokes2w} and the nonlinear Darcy problem~\eqref{eq:Darcy1w}-\eqref{eq:Darcy2w} with $\lambda=\lambda^{k}$:
			$$ \Psi (\lambda^{k}) = \iS_{f}(\lambda^{k}) + \iS_{p}(\lambda^{k}).
			$$
			\item Solve the linearized interface problem with a Krylov-type method (e.g., GMRES):
			$$ \int_{0}^{T} \left \langle \Psi^{\prime}(\lambda^{k})(h^{k}), \zeta \right \rangle_{\Gamma} = \int_{0}^{T} \left \langle -\Psi(\lambda^{k}), \zeta \right \rangle_{\Gamma}, \quad \forall \zeta \in \HL{\Lambda}.
			$$
			where the left-hand side is given by 
			$$ \Psi^{\prime}(\lambda^{k})(h^{k}) = S_{f,\lambda^{k}}^{\text{lin}}(h^k)+ S_{p, \lambda^{k}}^{\text{lin}}(h^k).$$
			That means each Krylov-iteration involves solution of linearized problems~\eqref{eq:linStokes1}-\eqref{eq:linDarcy2} to compute the matrix-free vector product on the left-hand side.  
			\item Update $\lambda^{k+1} = \lambda^{k} + h^{k}$, $k=k+1$, $\text{error}=\| h^{k}\|$.
		\end{algorithmic}

%The linearized interface problem~\eqref{eq:linIF} can be preconditioned by using the inverse operators of $S_{f,\lambda}^{\text{lin}}$ and $S_{p,\lambda}^{\text{lin}}$. That corresponds to solving the linearized Stokes problem with given normal velocity on the interface as Dirichlet boundary condition, and the linearized Darcy problem with given normal flux on the interface as Neumann boundary condition. 

The linearized interface problem~\eqref{eq:linIF} can be preconditioned by using the inverse operator of $S_{f,\lambda}^{\text{lin}}$ as proposed for the stationary case in \cite{Disca02}. That corresponds to solving the linearized Stokes problem with given normal velocity on the interface as Dirichlet boundary condition, and computing the normal stress on the interface. 

%
% ------------------------------------------------
%
%    SECTION: Nonconforming discretization in time
%	
% ------------------------------------------------
%	
\section{Nonconforming discretization in time} \label{sec:time}
As we solve the nonlinear interface problem~\eqref{eq:nonIF} globally in time, different time discretization schemes and/or different time step sizes can be used in the Stokes and Darcy regions. At the space-time interface, data is transferred from one space-time subdomain to a neighboring subdomain by using a suitable projection. 

\begin{figure}[ht!]
\centering
\includegraphics[scale=0.1]{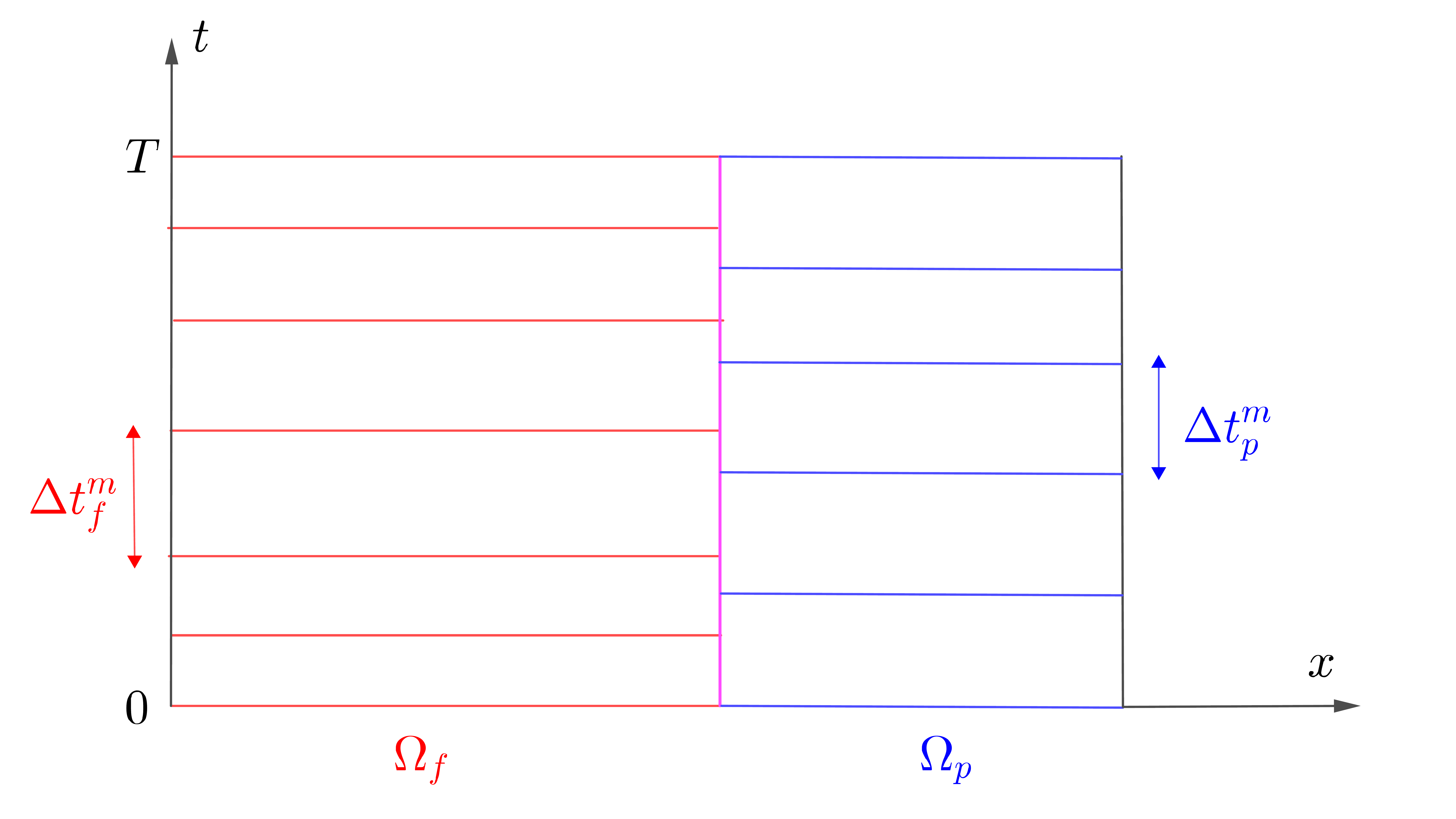} \vspace{-0.6cm}
\caption{Nonconforming time grids for the Stokes and Darcy problems.}
\label{fig:timegrids}
\end{figure}

We consider semi-discrete problems in time with nonconforming time grids. Let $ \iT_{f} $ and $ \iT_{p} $
be two possibly different partitions of the time interval $ (0,T) $ into sub-intervals (see Figure \ref{fig:timegrids}):
\begin{align*}
\iT_{f} &= \cup_{m=1}^{M} J_{f}^{m}, \; \text{with} \; J_{f}^{m}=(t_{f}^{m-1}, t_{f}^{m}], \quad \text{and} \;  \iT_{p} = \cup_{n=1}^{N} J_{p}^{n}, \; \text{with} \; J_{p}^{n}=(t_{f}^{n-1}, t_{f}^{n}].
\end{align*}
The time step sizes are $\Delta t_{f}^{m}= t^{m}-t^{m-1}$, $m=1, \hdots, M$, and $\Delta t_{p}^{n}= t^{n}-t^{n-1}$, $n=1, \hdots, N$, in the Stokes and Darcy regions, respectively. To simplify the discussion, the same temporal discretization scheme is considered for both subproblems; we use the backward Euler method for the time discretization and obtain the following semi-discrete local problems for the free flow
\begin{align}
&(\frac{\bu_{f}^{m}-\bu_{f}^{m-1}}{\Delta t_{f}^{m}}, \bv_{f})_{\Omega_{f}} + a_{f}(\bu_{f}^{m}, \bv_{f}) - b_{f}(\bv_{f},p_{f}^{m}) \nonumber\\
& \hspace{4.8cm} = (\bff_{f}^{m}, \bv_{f})_{\Omega_{f}} - b_{\Gamma f}(\bv_{f}, \lambda^{m}) , & \forall \bv_{f} \in \bX_{f}, \label{eq:Stokes1m}\\
& \hspace{3.3cm}  b_{f}(\bu_{f}^{m},q_{f}) = 0, & \forall q_{f} \in Q_{f}, \label{eq:Stokes2m}
\end{align}
and the Darcy flow
\begin{align}
a_{p}(\bu_{p}^{n}, \bv_{p}) - b_{p}(\bv_{p},p_{p}^{n}) &= - b_{\Gamma p}(\bv_{p}, \lambda^{n}), & \forall \bv_{p} \in \bX_{p}, \label{eq:Darcy1m}\\
b_{p}(\bu_{p}^{n}, q_{p})+\left (S_{p} \frac{p_{p}^{n}-p_{p}^{n-1}}{\Delta t_{p}^{n}}, q_{p}\right ) &= (f_{p}^{n}, q_{p})_{\Omega_{p}}, & \forall q_{p} \in Q_{p}. \label{eq:Darcy2m}
\end{align}
The wellposedness of the decoupled semi-discrete Stokes and Darcy problems \eqref{eq:Stokes1m} - \eqref{eq:Darcy2m} is followed from the strong monotonicity of the viscosity functions \eqref{eq:visf} and \eqref{eq:visp}, and \cite{Ervin14}.
The same idea can be generalized to higher order methods \cite{Japhet12}. 
	
For $i = f$ or $i=p$, we denote by $ P_{0}(\mathcal{T}_{i}, \HL{\Lambda}) $ the space of piecewise constant functions in time on grid $ \mathcal{T}_{i} $ with values in $  \HL{\Lambda} $:
\begin{equation} \label{P0space}
\begin{array}{ll}
P_{0}(\mathcal{T}_{f},  \HL{\Lambda}) &= \left \{ \phi: (0,T) \rightarrow  \HL{\Lambda},
\phi \text{ is constant on } J_{f}^{m}, \; \forall m=1, \dots, M \right \}, \vspace{2pt} \\
P_{0}(\mathcal{T}_{p},  \HL{\Lambda}) &= \left \{ \phi: (0,T) \rightarrow  \HL{\Lambda},
\phi \text{ is constant on } J_{p}^{n}, \; \forall n=1, \dots, N \right \}.
\end{array}
\end{equation}
In order to exchange data on the space-time interface between different time grids, we define the following
$ L^{2} $ projection $ \Pi_{p,f} $ from  $ P_{0} (\mathcal{T}_{f},  \HL{\Lambda}) $ onto $ P_{0}(\mathcal{T}_{p}, \HL{\Lambda}) $
(see \cite{OSWRwave,Japhet12})~: for $ \phi \in P_{0} (\mathcal{T}_{f},  \HL{\Lambda})$,
$ \Pi_{p,f} \phi \hspace{-2pt} \mid_{J_{p}^{n}} $ is the average value of $ \phi $ on $ J_{p}^{n} $,
for $ n=1, \dots, N $: \vspace{-0.2cm}
\begin{equation*}
\Pi_{p,f} \left ( \phi \right )\mid_{J_{p}^{n}}
  = \frac{1}{\mid J_{p}^{n}\mid} \sum_{m=1}^{M} \int_{J_{p}^{n} \cap J_{f}^{m}} \phi. 
\end{equation*}
The projection $ \Pi_{f,p} $ from  $ P_{0} (\mathcal{T}_{p},  \HL{\Lambda}) $ onto $ P_{0}(\mathcal{T}_{f}, \HL{\Lambda}) $ can be defined similarly. 
We use the algorithm described in \cite{Gander05} for effectively performing these projections.

Next, we weakly enforce the transmission conditions over the time intervals with nonconforming time grids. We still denote by $ (\bu_{f}, p_{f})$ and $(\bu_{p}, p_{p})$ the solution of the semi-discrete in time problems. We choose $ \lambda $ piecewise constant in time on one grid, either $ \iT_{f} $ or $ \iT_{p} $. For the Stokes-Darcy coupling, the flow is supposed to be faster in the fluid domain than that in the porous medium, thus we choose $\lambda \in P_{0} (\iT_{f},  \HL{\Lambda})$ and impose 
$$ \left ( -\bn_{f} \cdot (\nu_{f}(\vert \bD(\bu_{f})\vert) \bD(\bu_{f})) \cdot \bn_{f} + p_{f}\right )\vert_{\Gamma}= \Pi_{f,f} (\lambda) = \lambda. $$ The weak continuity
of the normal stress in time across the interface
is fulfilled by letting 
$$
p_{p}\vert_{\Gamma} = \Pi_{p,f} (\lambda ) \in P_{0} (\iT_{p}, \HL{\Lambda}).
$$
The semi-discrete (nonconforming in time) counterpart of the normal velocity continuity~\eqref{eq:TC1}
is weakly enforced by integrating it over each time interval $ J_{f}^{m} $ of grid $ \iT_{f} $ :
$\forall m=1,...,M$,
\begin{equation} \label{eq:nonIFtime}
\int_{J^{m}_{f}}  \biggl( \left \langle S_{f}(\lambda), \zeta \right \rangle_{\Gamma} + \left \langle \Pi_{f,p} \bigl ( S_{p} \left (\Pi_{p,f}(\lambda)\right ) \bigl), \zeta \right \rangle_{\Gamma} \biggl ) \, ds= 0, \quad \forall \zeta \in  \HL{\Lambda}.
\end{equation}
Similarly for the linearized interface problem:
\begin{align}
&\int_{J^{m}_{f}} \biggl( \left \langle S_{f,\lambda^{k}}^{\text{lin}}(h^{k}), \zeta \right \rangle_{\Gamma}+ \left \langle \Pi_{f,p} \bigl ( S_{p, \lambda^{k}}^{\text{lin}} (\Pi_{p,f}(h^{k}) \bigl), \zeta \right \rangle_{\Gamma} \biggl) \, ds = \nonumber\\
& \hspace{1cm}  \int_{J^{m}_{f}} \biggl( \left \langle -S_{f}(\lambda^{k}), \zeta \right \rangle_{\Gamma} + \left \langle -\Pi_{f,p} \bigl ( S_{p} \left (\Pi_{p,f}(\lambda^{k})\right ) \bigl), \zeta \right \rangle_{\Gamma} \biggl ) \, ds, \quad \forall \zeta \in  \HL{\Lambda}. \label{eq:linIFtime}
\end{align}

% ------------------------------------------------
%
%    SECTION 5 : Numerical results
%	
% ------------------------------------------------
%	
	
\section{Numerical results}
\label{sec:NumRe}
We investigate the numerical performance of the proposed global-in-time decoupling algorithm on two test cases: Test case 1 with a known solution and Test case 2 where the flow is driven by a pressure drop. For the latter, we consider both continuous and discontinuous parameters. We shall verify the accuracy in space and in time, and the efficiency of the proposed method with nonconforming time grids over conforming time grids. Note that the code to generate the results below is implemented in {\em FreeFem++}~\cite{Freefem} in a sequential setting, and we do not investigate parallel performance of the method in this work.
%
%
%  TEST CASE 1 (Hyesuk's paper)
%
%
\subsection{Test case 1: with a known analytical solution}
We consider a test case with a known exact solution. The fluid domain and porous medium are $\Omega_{f}=(0,1) \times (1,2)$ and $\Omega_{p}=(0,1) \times (0,1)$ respectively, and the exact solution is given by
\begin{align*}
\bu_{f}&=\left [ (y-1)^{2}x^{3}(1+t^{2}), \; -\cos(y) e(1+t^{2}) \right ], \\
p_{f}&=\left (\cos(y)e^{y}+y^{2}-2y+1\right ) \left (1+t^{2}\right ), \\
\bu_{p} & = \left [-x \left (\sin(y)e + 2(y-1)\right ) \left (1+t^{2}\right ), \; \left (-\cos(y) e + (y-1)^{2}\right ) \left (1+t^{2}\right ) \right ], \\
p_{p}&=\left (-\sin(y) e+\cos(x) e^{y}+y^{2}-2y+1\right ) \left (1+t^{2}\right ),
\end{align*}
for which the Beavers-Joseph-Saffman condition is satisfied with $\alpha=1$. We perform the numerical experiments with the following parameters: $\kappa=1$, $K_{f}=K_{p}=1$, $\nu_{f\infty}=\nu_{p\infty}=~0.5$ and $\nu_{f0}=\nu_{p0}=1.5$. 
The boundary and initial conditions are imposed using the exact solution. For finite element approximations, we consider structured meshes and use either (i) the Taylor-Hood elements for both Stokes and Darcy problems or (ii) the MINI elements for the Stokes problem and the Raviart-Thomas of order~1 elements for the Darcy problem. In addition, a stability term $\eta \left (\nabla \cdot \bu_{p}, \nabla \cdot \bv_{p}\right )$ was added to the Darcy equation with $\eta=10$ as the exact Darcy velocity field is divergence free. 

We shall verify the convergence rates in space and in time of the proposed algorithm with nonconforming time grids.
For the iterative solvers, unless otherwise specified, only one Newton iteration is performed (i.e., k=1 in {Algorithm 1}) and GMRES stops when the relative residual is smaller than the tolerance $\varepsilon=10^{-7}$ or when the maximum number of iterations, itermax$=100$, is reached. 
%
%%%%%%%%%%%%%%%%%%%%%%%%%%%%%%%%%%%%%%%%%%%
%
% CONVERGENCE IN Space
%
%%%%%%%%%%%%%%%%%%%%%%%%%%%%%%%%%%%%%%%%%%%
%
We first investigate the accuracy in space for both linear viscosities with $r_{f}=r_{p}=2$ and nonlinear viscosities with $r_{f}=r_{p}=1.5$. Tables~\ref{tab:THConvSpace} and \ref{tab:RTConvSpace} show the errors at $T=0.01$ with $\Delta t_{f}=0.002$ and $\Delta t_{p}=0.001$ for the linear and nonlinear problems using different finite element spaces. As this is a non-physical example, we have chosen a large time step in the fluid domain and a small time step in the porous medium. In the next test case, we will consider the choice where the time step size in the fluid domain is smaller. We observe from Tables~\ref{tab:THConvSpace} and \ref{tab:RTConvSpace}  that the orders of accuracy in space are preserved with nonconforming time grids. 
In addition, concerning the convergence of GMRES to solve the linearized interface problem, we show in Table~\ref{tab:RTConvSpaceiter} the number of GMRES iterations needed to reach the tolerance $\varepsilon=10^{-10}$ for the case with no preconditioner and with the preconditioner $\left (S_{f,\lambda}^{\text{lin}}\right )^{-1}$. First, we notice the number of iterations required is reasonable; for the case without preconditioner, it is increasing slightly when the mesh size is decreasing while for the preconditioned system, the number of iterations remain small when $h$ is small.

\begin{table}[http!]
	\small
	\setlength{\extrarowheight}{4pt}
	\centering
	\begin{tabular}{| l | l | l l | l l | l l | l l |} \hline 
		\multicolumn{2}{|c|}{$h$} & \multicolumn{2}{c|}{$1/4$} & \multicolumn{2}{c|}{$1/8$} & \multicolumn{2}{c|}{$1/16$} & \multicolumn{2}{c|}{$1/32$}\\ \hline
		\multicolumn{10}{|c|}{Linear viscosities}  \\ \hline
		%
		%\multicolumn{2}{|c|}{Number of iterations}              & \multicolumn{2}{c|}{100} & \multicolumn{2}{c|}{100} & \multicolumn{2}{c|}{100} & \multicolumn{2}{c|}{100}  \\ \hline
		%
		\multirow{2}{*}{$\bu_{f}$} & $L^{2}$ error  			& 9.07e-04    &    & 9.33e-05   & [3.28]  &  1.19e-05& [2.97]  & 1.79e-06 & [2.73]    \\ 
												& $H^{1}$ error 			&	2.64e-02	 &    &  5.55e-03  & [2.25]  & 1.38e-03 & [2.01]  & 3.72e-04  & [1.89]   \\  \hline
		$p_{f}$ 								&  $L^{2}$ error  			&	2.91e-02	 &    &  5.55e-03  & [2.39]  & 1.36e-03 & [2.03]   &  3.97e-04 & [1.78]   \\ \hline
		\multirow{2}{*}{$\bu_{p}$} & $L^{2}$ error  			  & 1.29e-03  &    &  1.71e-04  & [2.92]   & 2.02e-05 & [3.08]  &  4.58e-06  & [2.14]  \\ 
												&  $H^{\text{div}}$  error & 2.24e-03	 &    &  3.42e-04  & [2.71]   & 8.19e-05  & [2.06]  & 1.93e-05 & [2.09]    \\ \hline
		$p_{p}$								&  $L^{2}$ error  			  & 3.12e-02	 &    &  5.07e-03  & [2.62]  & 1.34e-03  & [1.92]  & 3.24e-04 & [2.05]    \\ \hline
		\multicolumn{10}{|c|}{Nonlinear viscosities} \\ \hline
		%
		%\multicolumn{2}{|c|}{Number of iterations}              & \multicolumn{2}{c|}{100} & \multicolumn{2}{c|}{100} & \multicolumn{2}{c|}{100} & \multicolumn{2}{c|}{100}  \\ \hline
		\multirow{2}{*}{$\bu_{f}$} & $L^{2}$ error  			& 9.64e-04    &    & 1.03e-04   & [3.23]  &  1.82e-05& [2.50]  & 1.26e-05 &    \\ 
												& $H^{1}$ error 			&	2.77e-02	 &    &  5.88e-03  & [2.24]  & 1.47e-03 & [2.00]  & 4.18e-04 &   [1.81]] \\ \hline
		$p_{f}$ 								&  $L^{2}$ error  			&	2.84e-02	 &    &  5.50e-03  & [2.37]   & 1.50e-03 & [1.88]   & 7.36e-04 & [1.03]   \\ \hline
		\multirow{2}{*}{$\bu_{p}$} & $L^{2}$ error  			  & 1.26e-03  &    &  1.72e-04  & [2.87]  & 2.10e-05 &  [3.03]  & 7.44e-06  &     \\ 
												&  $H^{\text{div}}$  error & 2.18e-03	 &    &  3.37e-04  & [2.69]  & 8.00e-05  & [2.08]  & 1.98e-05  &   [2.02] \\ \hline
		$p_{p}$								&  $L^{2}$ error  			  & 3.02e-02	 &    &  4.91e-03  & [2.62]  & 1.30e-03  & [1.92]  & 3.15e-04 &  [2.05]  \\ 	\hline							
	\end{tabular}  
\caption{[Test case 1] Errors with Taylor-Hood elements for the Stokes and Darcy problems at $T=0.01$ with $\Delta t_{f}=0.002$ and $\Delta t_{p}=0.001$.}  \label{tab:THConvSpace}
\end{table}

\begin{table}[http!]
	\small
	\setlength{\extrarowheight}{4pt}
	\centering
	\begin{tabular}{| l | l | l l | l l | l l | l l |} \hline 
		\multicolumn{2}{|c|}{$h$} & \multicolumn{2}{c|}{$1/4$} & \multicolumn{2}{c|}{$1/8$} & \multicolumn{2}{c|}{$1/16$} & \multicolumn{2}{c|}{$1/32$}\\ \hline
		\multicolumn{10}{|c|}{Linear viscosities}  \\ \hline
		%
		%\multicolumn{2}{|c|}{Number of iterations}              & \multicolumn{2}{c|}{12} & \multicolumn{2}{c|}{17} & \multicolumn{2}{c|}{23} & \multicolumn{2}{c|}{34}  \\ \hline
		\multirow{2}{*}{$\bu_{f}$} & $L^{2}$ error  			& 1.09e-02    &    &  2.63e-03  & [2.05]  & 6.72e-04 & [1.97]  & 1.82e-04 & [1.89]   \\ 
												& $H^{1}$ error 			&	2.24e-01	&    &   9.88e-02 & [1.18]   & 5.02e-02 & [0.98]  & 2.69e-02 & [0.90]     \\  \hline
		$p_{f}$ 								&  $L^{2}$ error  			&	2.51e-01	 &    &  6.65e-02  & [1.92]   & 2.22e-02& [1.58]  & 1.06e-02 & [1.07]   \\  \hline
		\multirow{2}{*}{$\bu_{p}$} & $L^{2}$ error  			  & 2.11e-02  &    &  4.29e-03 &[2.30]   &  1.10e-03& [1.96]  &  2.65e-04 & [2.05]   \\  
												&  $H^{\text{div}}$  error & 2.36e-02	&    &  5.03e-03  & [2.23]  & 1.28e-03 & [1.97]  & 3.12e-04 & [2.04]   \\  \hline
		$p_{p}$								&  $L^{2}$ error  			  & 3.04e-02	&    &  4.94e-03  & [2.62]  & 1.31e-03  & [1.92]  & 3.16e-04  &  [2.05]    \\  \hline
		\multicolumn{10}{|c|}{Nonlinear viscosities} \\ \hline
		%
		%\multicolumn{2}{|c|}{Number of iterations}              & \multicolumn{2}{c|}{12} & \multicolumn{2}{c|}{16} & \multicolumn{2}{c|}{21} & \multicolumn{2}{c|}{30}  \\ \hline
		\multirow{2}{*}{$\bu_{f}$} & $L^{2}$ error  			& 1.09e-02    &    &  2.62e-03  & [2.05]  & 6.70e-04 &  [1.97] & 1.81e-04 & [1.89]   \\ 
												& $H^{1}$ error 			&	2.24e-01	&    &   9.88e-02 &  [1.18]  & 5.02e-02 &  [0.98] & 2.69e-02 & [0.90]    \\  \hline
		$p_{f}$ 								&  $L^{2}$ error  			&	2.07e-01	 &    &  5.40e-02  & [1.94]  & 1.86e-02& [1.54]  & 8.79e-03 & [1.08]   \\  \hline
		\multirow{2}{*}{$\bu_{p}$} & $L^{2}$ error  			  & 2.12e-02  &    &  4.30e-03 & [2.30]  &  1.10e-03& [1.97]   &2.66e-04  & [2.05]    \\  
												&  $H^{\text{div}}$  error & 2.37e-02	&    &  5.02e-03  & [2.24]  & 1.28e-03 & [1.97]   & 3.12e-04 & [2.04   \\  \hline
		$p_{p}$								&  $L^{2}$ error  			  & 2.95e-02	&    &  4.79e-03  & [2.63]  & 1.27e-03  & [1.92] & 3.08e-04  & [2.04]     \\  \hline							
	\end{tabular}  
\caption{[Test case 1] Errors with MINI elements for the Stokes problem and with Raviart-Thomas of order 1 elements for the Darcy problem at $T=0.01$ with $\Delta t_{f}=0.002$ and $\Delta t_{p}=0.001$.}  \label{tab:RTConvSpace}
\end{table}

\begin{table}[http!]
	\small
	\setlength{\extrarowheight}{4pt}
	\centering
	\begin{tabular}{| c |  p{0.7cm} |   p{0.7cm} |   p{0.7cm}  |   p{0.7cm} | p{0.7cm} |   p{0.7cm} |   p{0.7cm}  |   p{0.7cm} |} \hline 
		  \multirow{2}{*}{$h$} &\multicolumn{4}{c|}{Linear viscosities}  & \multicolumn{4}{c|}{Nonlinear viscosities} \\ \cline{2-9}
		& $1/4$ & $1/8$ & $1/16$ & $1/32$ & $1/4$ & $1/8$ & $1/16$ & $1/32$\\ \hline
		With no preconditioner       & 17 & 24 & 32 & 46  & 16 & 23 & 30 & 44  \\ \hline
		With a preconditioner         & 21 & 22 & 17 & 21 & 23 & 25 & 18 & 18  \\ \hline
		\end{tabular}  
\caption{[Test case 1] Number of GMRES iterations needed to reach the tolerance $10^{-10}$ using MINI elements for the Stokes problem and with Raviart-Thomas of order 1 elements for the Darcy problem at $T=0.01$ with $\Delta t_{f}=0.002$ and $\Delta t_{p}=0.001$.}  \label{tab:RTConvSpaceiter}
\end{table}
%%%%%%%%%%%%%%%%%%%%%%%%%%%%%%%%%%%%%%%%%%%
%
% CONVERGENCE IN TIME
%
%%%%%%%%%%%%%%%%%%%%%%%%%%%%%%%%%%%%%%%%%%%

For time errors, we analyze the accuracy in time when nonconforming time grids are used. Toward this end, we fix $h=1/32$ and denote by $ \Delta t_{\text{coarse}} \in \{ 0.2, 0.1, 0.05, 0.025 \}$ the coarse time step sizes, and $\Delta t_{\text{fine}}=\Delta t_{\text{coarse}}/2$ the fine time step size.  We consider three types of time grids as follows:
\begin{itemize}
\item[i)] Coarse conforming time grids: $\Delta t_{f} = \Delta t_{p} =\Delta t_{\text{coarse}}$.
\item[ii)] Fine conforming time grids: $\Delta t_{f} = \Delta t_{p} = \Delta t_{\text{fine}}$.
\item[iii)] Nonconforming time grids: $\Delta t_{f}=\Delta t_{\text{coarse}}$ and $\Delta t_{p}=\Delta t_{\text{fine}}$.
\end{itemize}

\begin{figure}[http!]
    \centering
    \begin{subfigure}[b]{0.45\textwidth}
        \includegraphics[width=\textwidth]{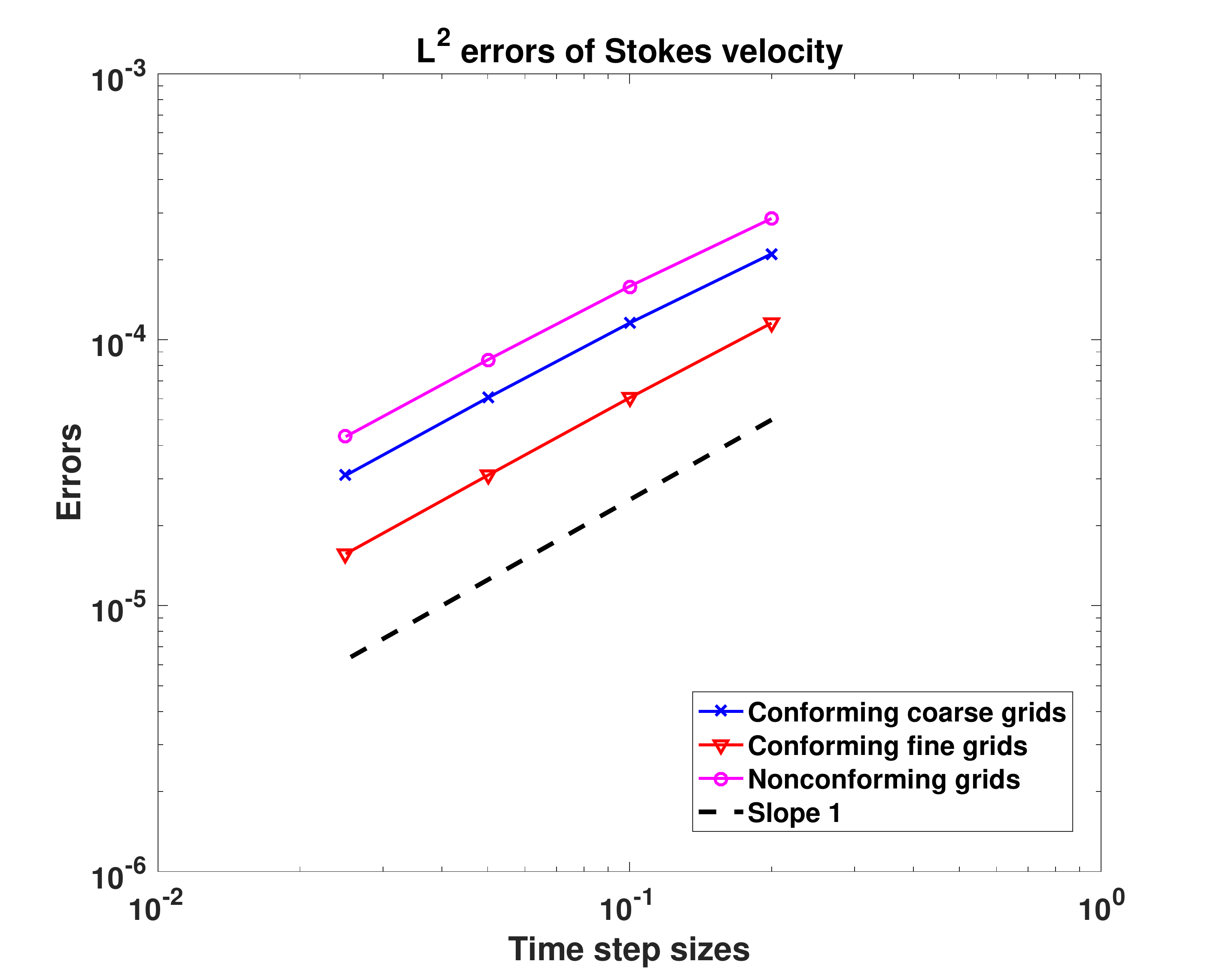}
    \end{subfigure}
    \begin{subfigure}[b]{0.45\textwidth}
        \includegraphics[width=\textwidth]{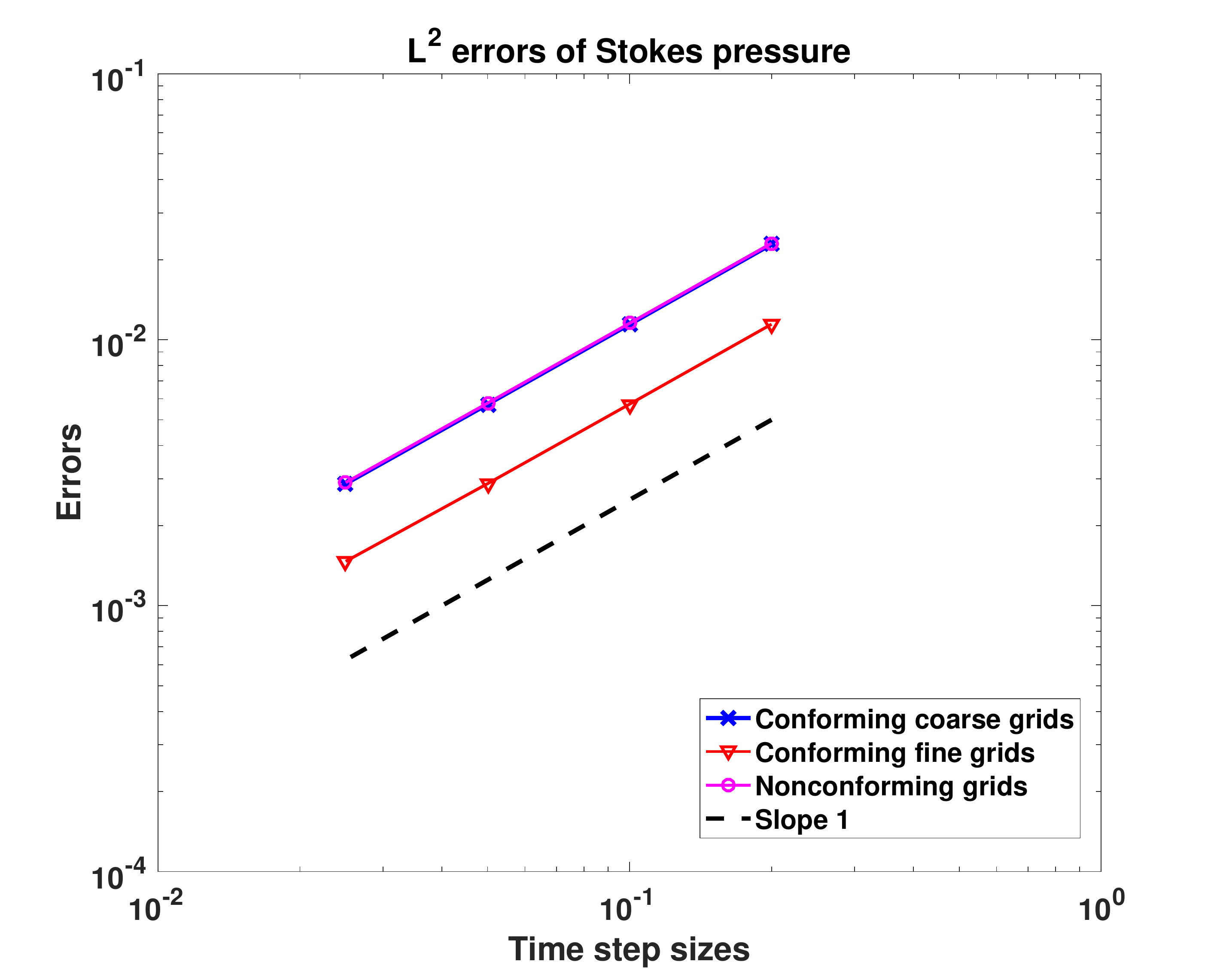}
    \end{subfigure} \vspace{0.3cm}\\
    \begin{subfigure}[b]{0.45\textwidth}
        \includegraphics[width=\textwidth]{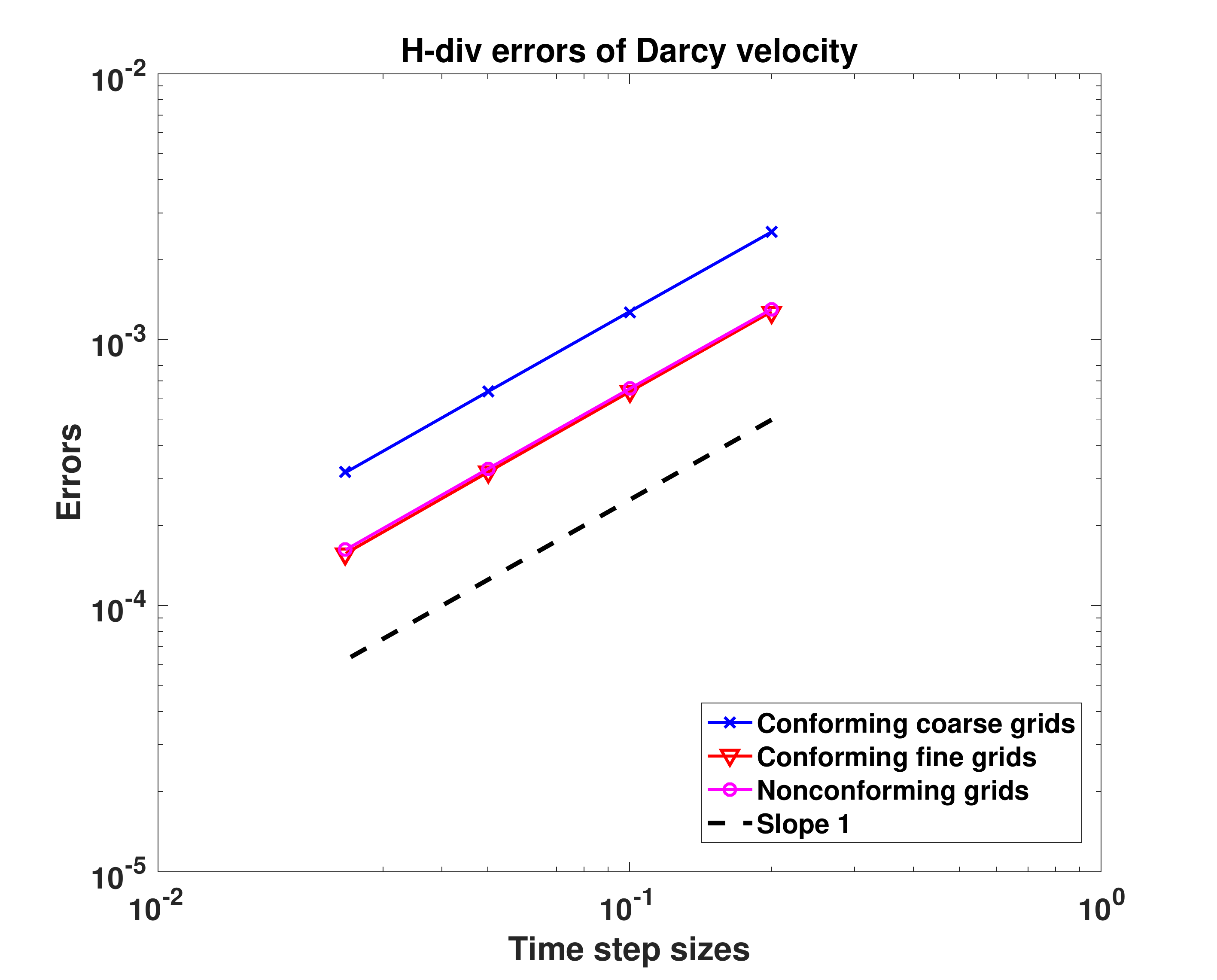}
    \end{subfigure}
    \begin{subfigure}[b]{0.45\textwidth}
        \includegraphics[width=\textwidth]{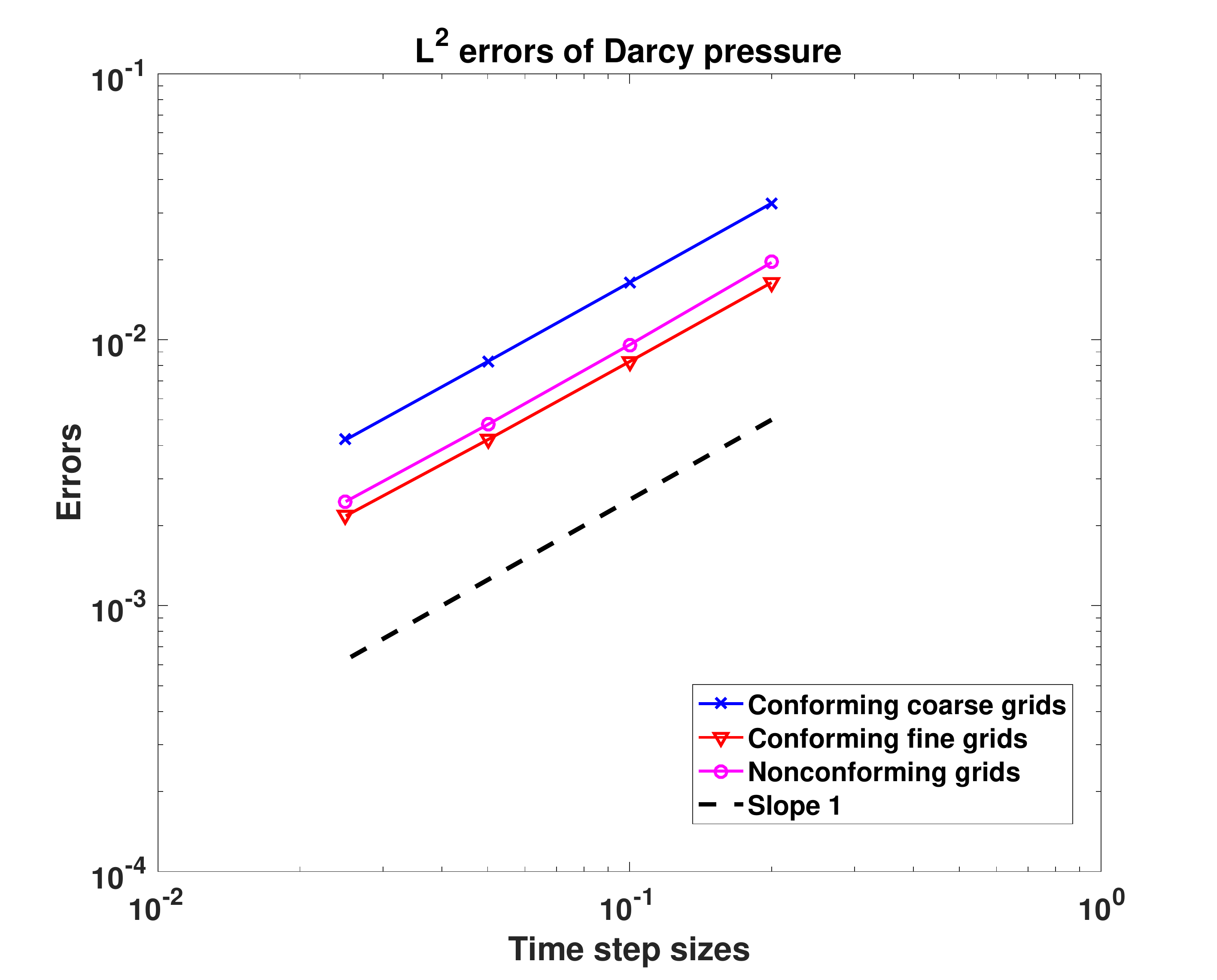}
    \end{subfigure} 
    \caption{[Test case 1] Errors for the linear Stokes and Darcy problems at $T=0.2$ with Taylor-Hood elements.}\label{fig:LinearConvTime}
\end{figure}

\begin{figure}[http!]
    \centering
    \begin{subfigure}[b]{0.45\textwidth}
        \includegraphics[width=\textwidth]{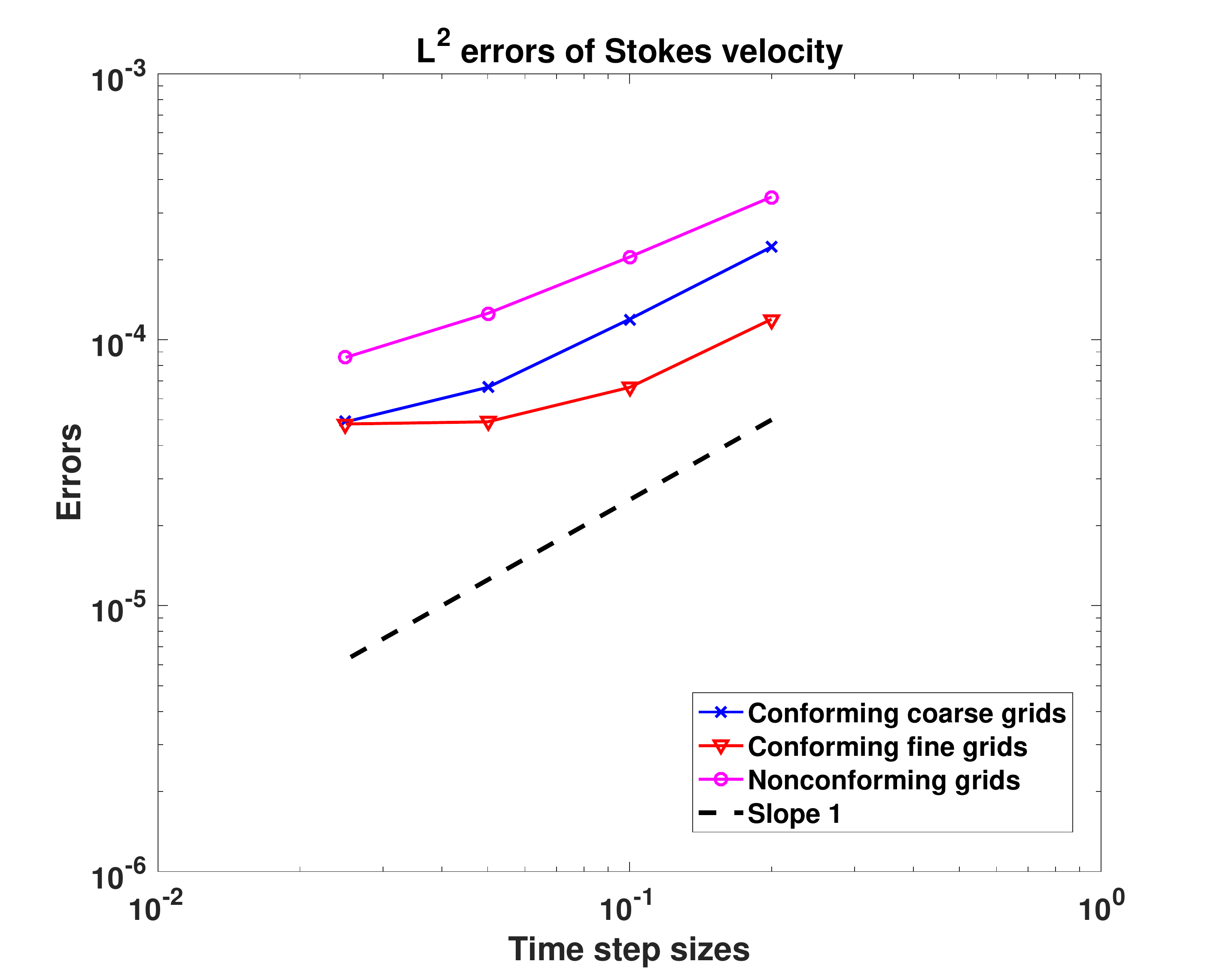}
    \end{subfigure}
    \begin{subfigure}[b]{0.45\textwidth}
        \includegraphics[width=\textwidth]{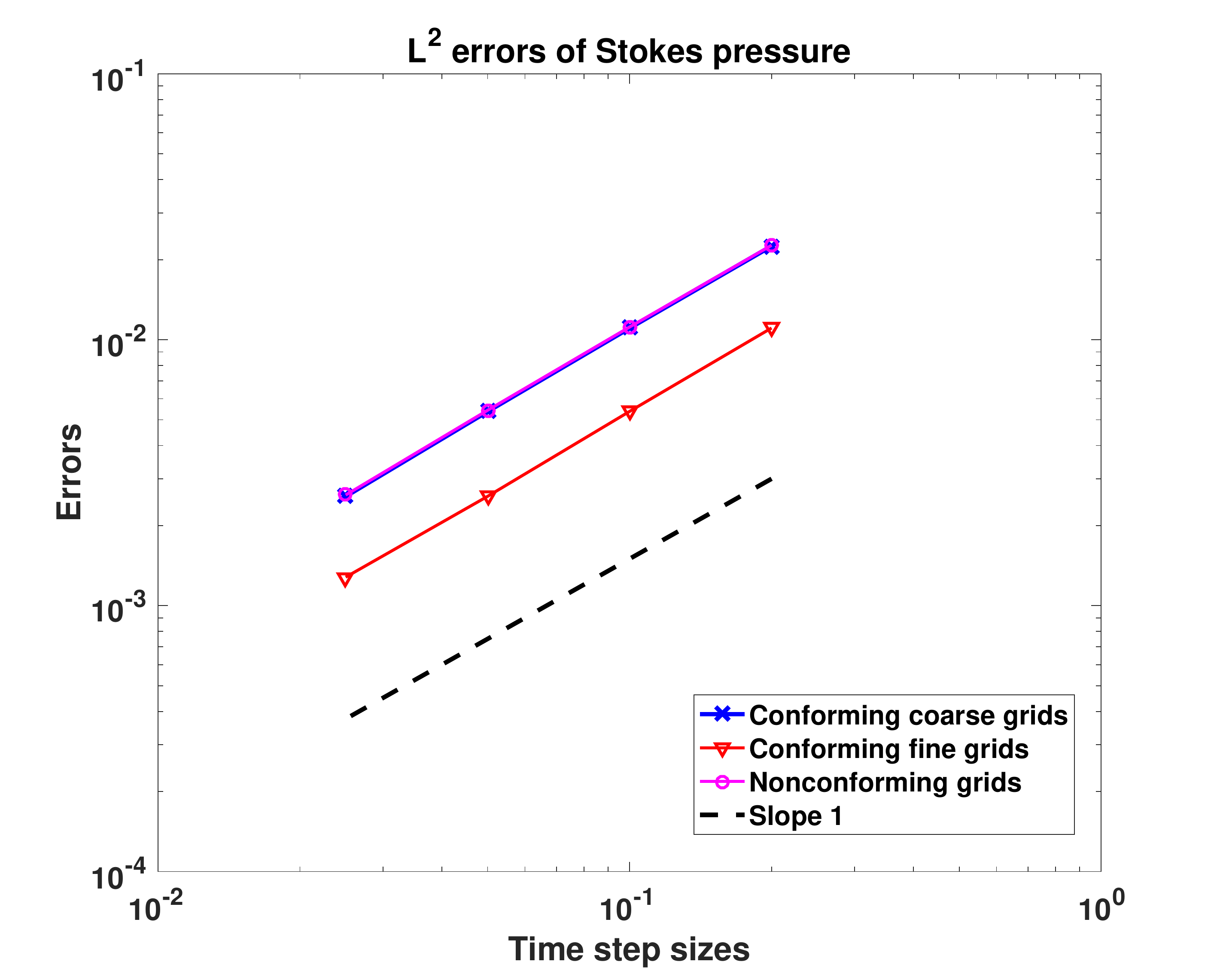}
    \end{subfigure} \vspace{0.3cm}\\
    \begin{subfigure}[b]{0.45\textwidth}
        \includegraphics[width=\textwidth]{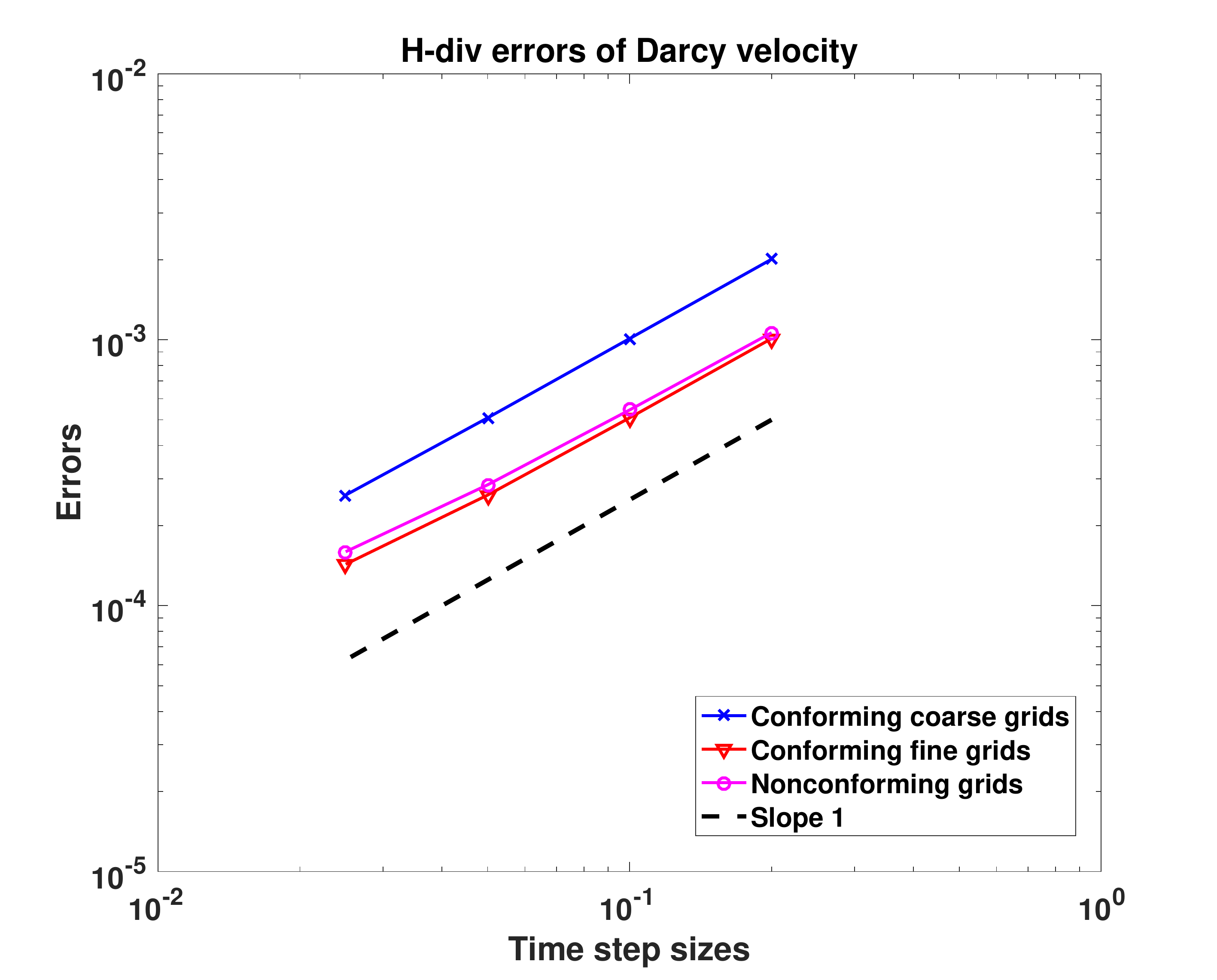}
    \end{subfigure}
    \begin{subfigure}[b]{0.45\textwidth}
        \includegraphics[width=\textwidth]{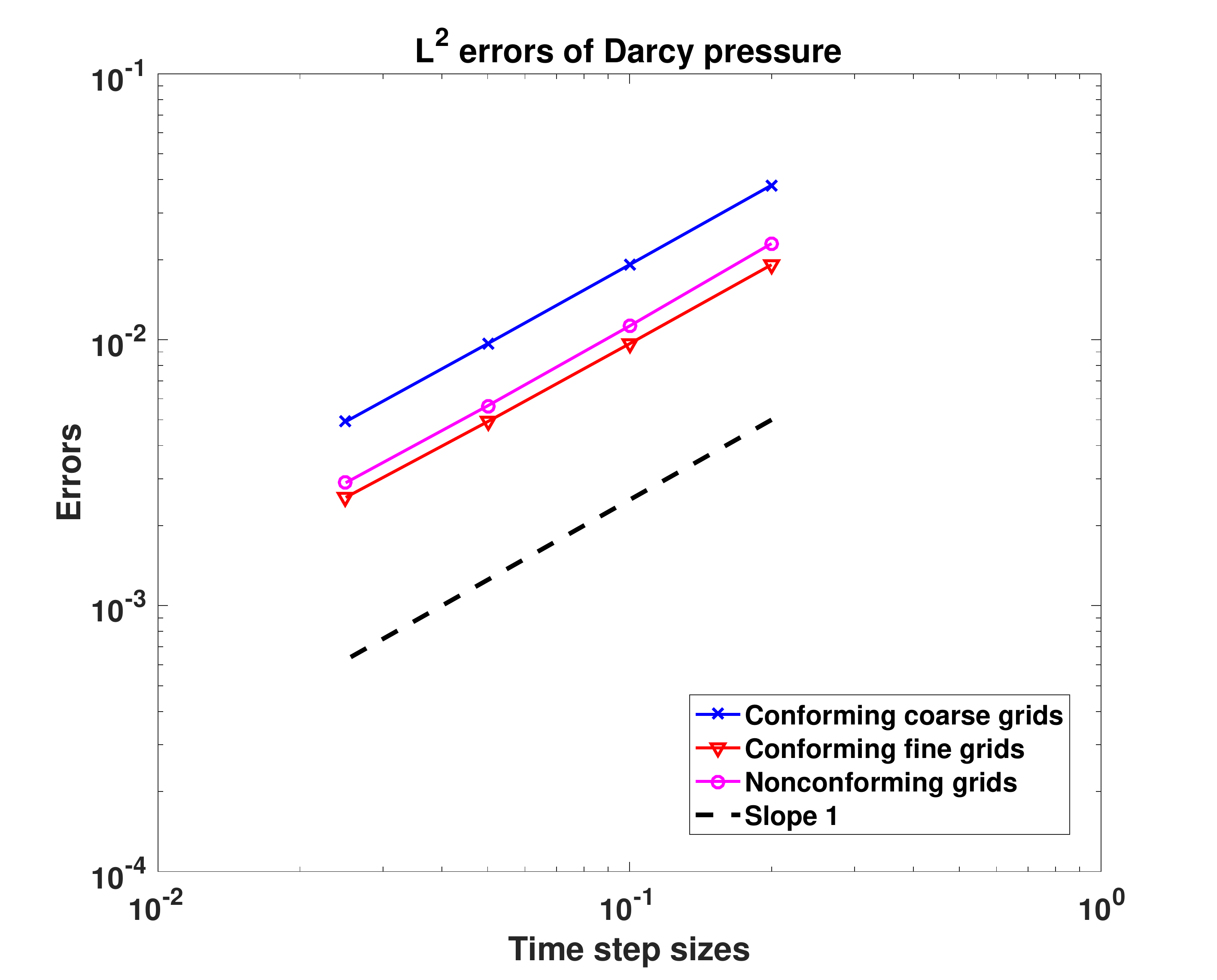}
    \end{subfigure} 
    \caption{[Test case 1] Errors for the nonlinear Stokes and Darcy problems at $T=0.2$ with Taylor-Hood elements.}\label{fig:nonLinearConvTime}
\end{figure}

We first consider the approximations by Taylor-Hood elements. Figures~\ref{fig:LinearConvTime} and \ref{fig:nonLinearConvTime} show the errors for the linear and nonlinear viscosities respectively. We observe
that first order convergence is preserved with the nonconforming time grids. Moreover,
the errors with nonconforming time grids (in magenta) in the porous medium are close to those with fine conforming time steps (in red), which is expected as a smaller time step is used in the porous medium. Likewise, the errors with nonconforming time grids (in magenta) in the fluid domain are close to those with coarse conforming time steps (in blue). Thus the accuracy in time of the solution is preserved with the nonconforming time grids. Moreover, in Table~\ref{tab:THcputime}, we compare the computer running times when using conforming and nonconforming time grids, which shows that using nonconforming time grids could significantly reduce the computational time while still maintaining the desired accuracy.

\begin{table}[http!]
\begin{center}
\begin{tabular}{| c ||  c | c | c | c |}
\firsthline
 \multirow{2}{*}{$\Delta t$}  &  \multicolumn{2}{c|}{Linear viscosities}  & \multicolumn{2}{c|}{Nonlinear viscosities}   \\[1pt] \cline{2-5}
    & Conforming & Nonconforming & Conforming & Nonconforming \\[1pt]
\hline
$0.2$  		&  87 	&   & 122  &      \\[1pt]   & & 143	&  & 178  \\[1pt]  
$0.1 $ 		&  209  &   & 287  &    \\[1pt]     & & 285	& 	&	348  \\[1pt]  
$0.05 $  	&  432 	&  & 578   &    \\[1pt]    & & 578	& 	&	710 \\[1pt]  
$0.025 $ 	&  893	&  & 1127  &   \\[1pt]     & & 1176	& 	&	1424  \\[1pt]  
$0.0125 $ & 1816  &  & 2175 &   \\
\lasthline
\end{tabular}
\end{center}
\caption{Comparison of the computer running times (in seconds) of conforming and nonconforming time grids with Taylor-Hood elements on a fixed mesh $h=1/32$. } \label{tab:THcputime}
\end{table}

\begin{figure}[http!]
    \centering
    \begin{subfigure}[b]{0.45\textwidth}
        \includegraphics[width=\textwidth]{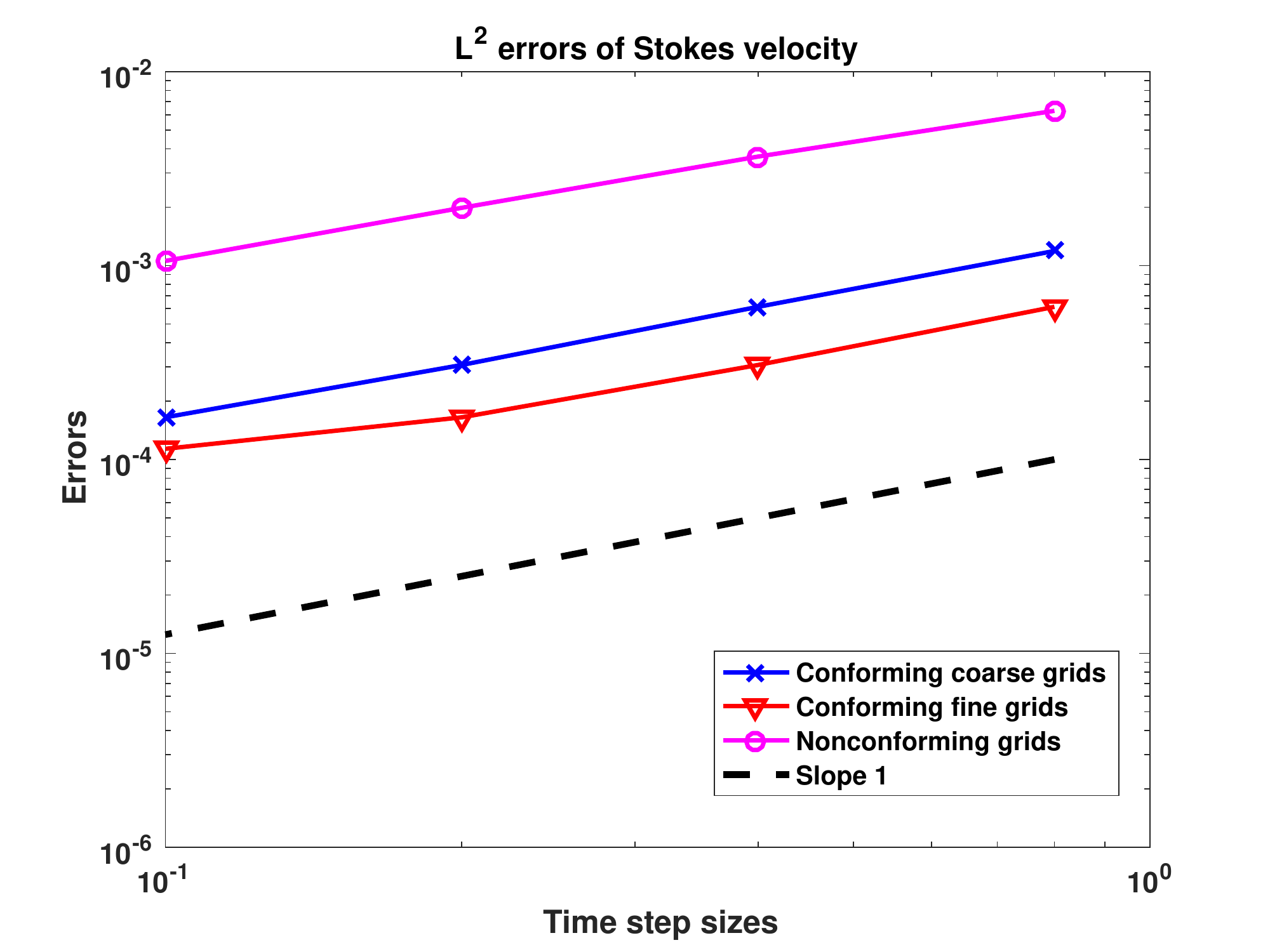}
    \end{subfigure}
    \begin{subfigure}[b]{0.45\textwidth}
        \includegraphics[width=\textwidth]{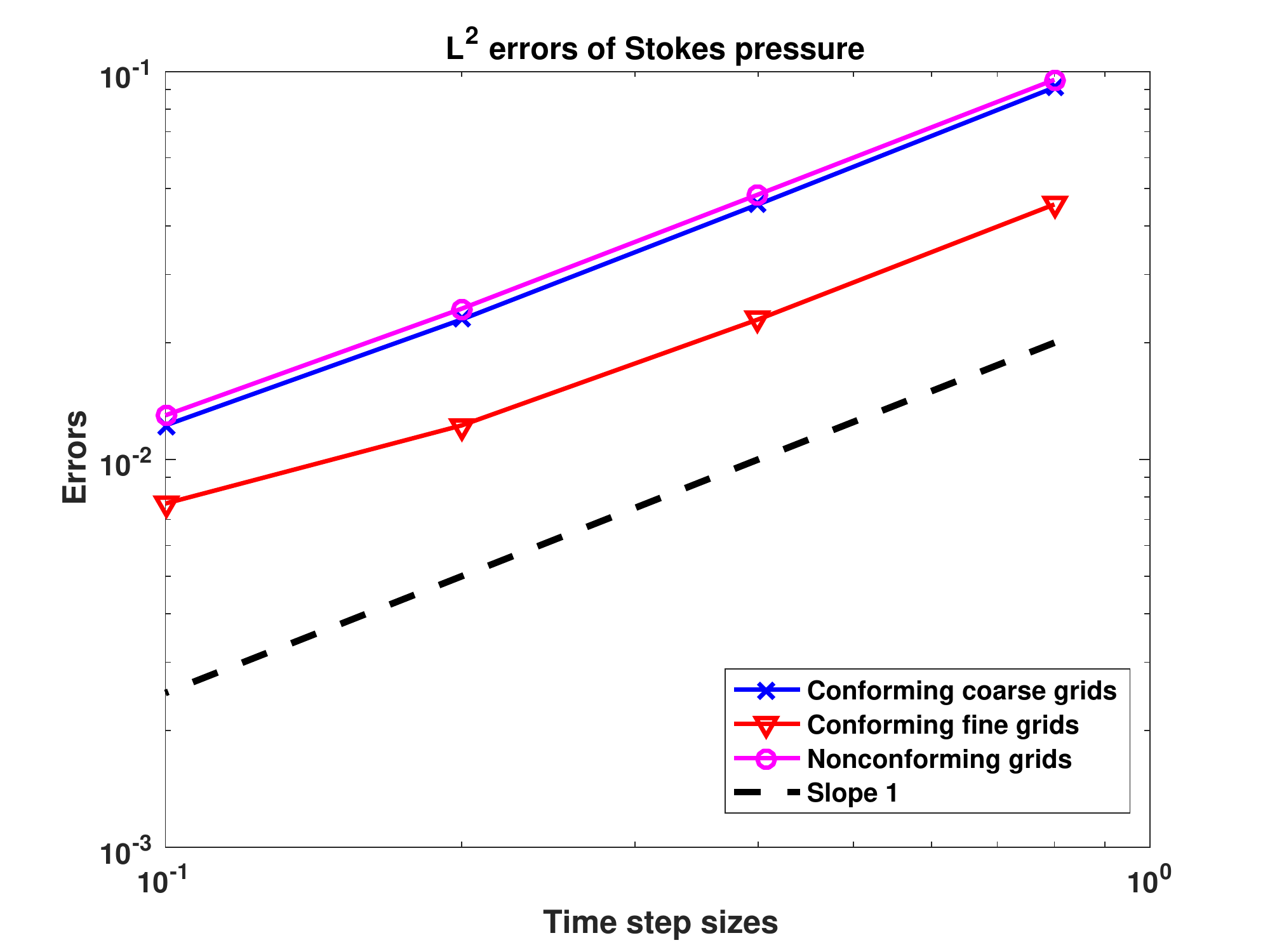}
    \end{subfigure} \vspace{0.3cm}\\
    \begin{subfigure}[b]{0.45\textwidth}
        \includegraphics[width=\textwidth]{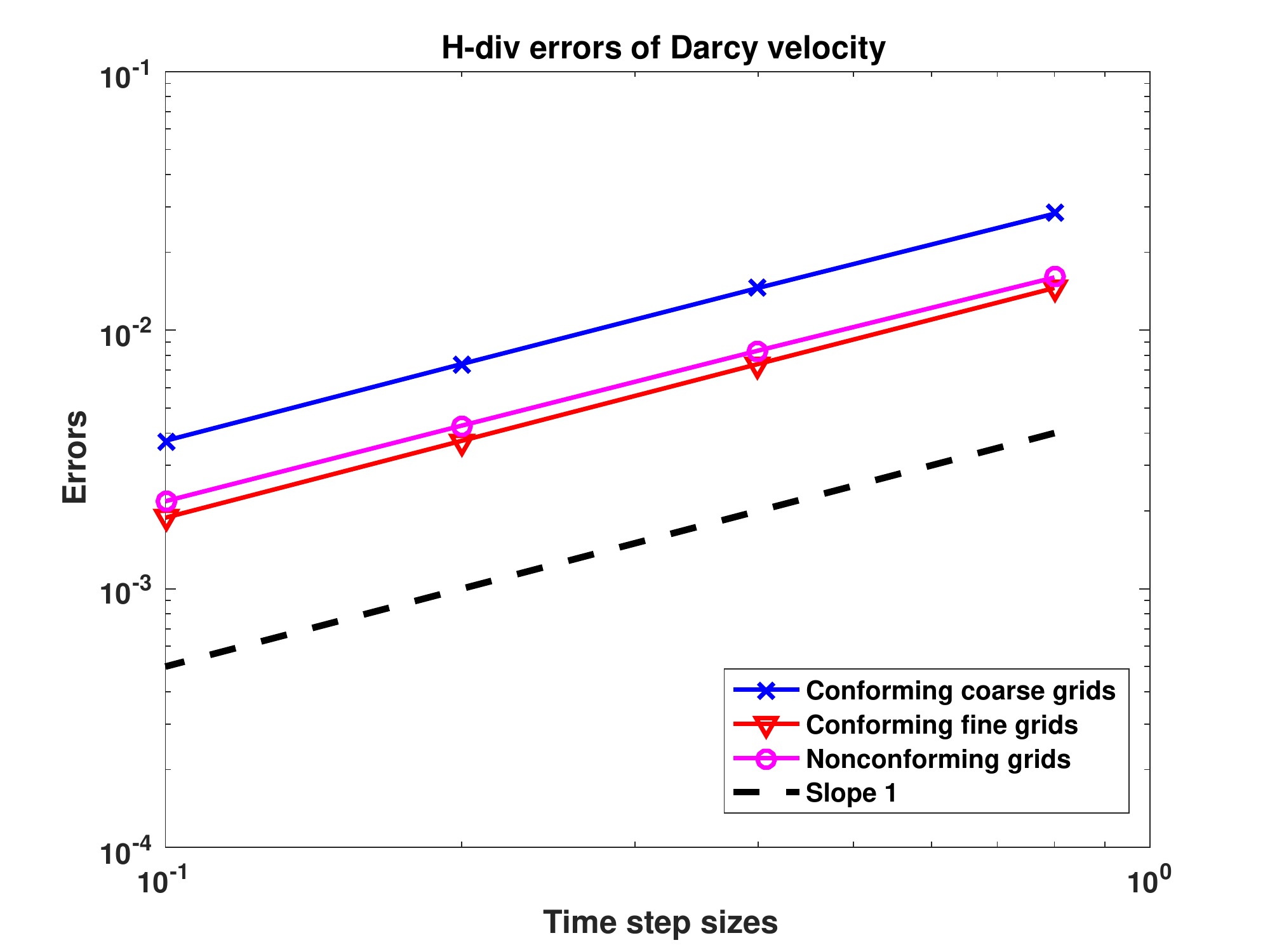}
    \end{subfigure}
    \begin{subfigure}[b]{0.45\textwidth}
        \includegraphics[width=\textwidth]{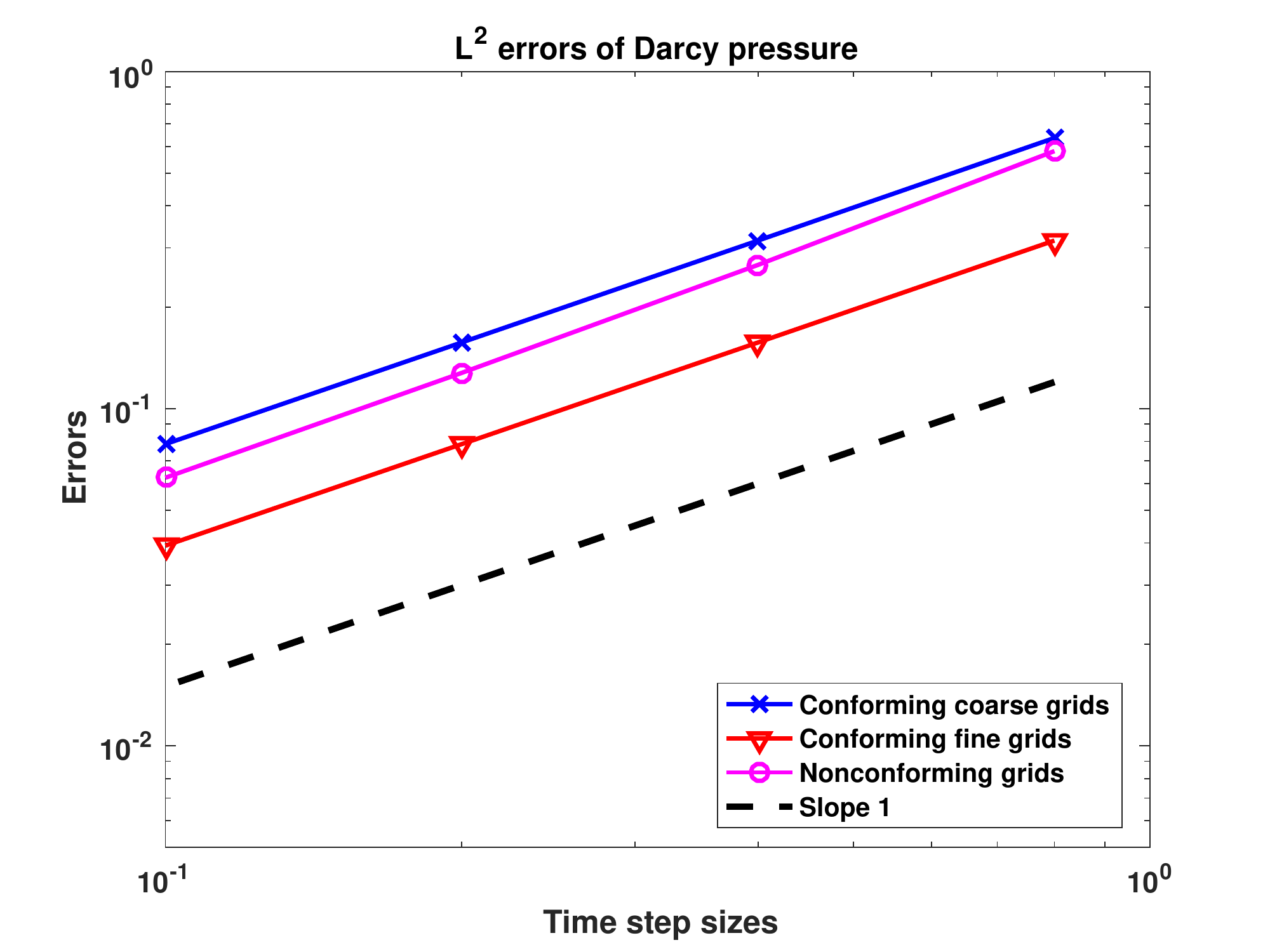}
    \end{subfigure} 
    \caption{[Test case 1] Errors for the nonlinear Stokes and Darcy problems at $T=0.8$ with MINI elements for the Stokes problem and Raviart-Thomas elements for the Darcy problem.}\label{fig:nonLinearConvTimeRT}
\end{figure}

We perform a similar test using MINI elements for the Stokes problem and Raviart-Thomas elements for the Darcy problem. We fix $h=1/64$, and consider $ \Delta t_{\text{coarse}} \in \{ 0.8, 0.4, 0.2, 0.1 \}$ and $\Delta t_{\text{fine}}=\Delta t_{\text{coarse}}/2$. The final time is large, $T=0.8$, thus we use two Newton iterations for the nonlinear solvers (instead of only one iteration). Figure~\ref{fig:nonLinearConvTimeRT} shows the errors for the case with nonlinear viscosities, which again confirms that the convergence order and accuracy in time are preserved with nonconforming time grids. In addition, we report in Table~\ref{tab:RTcputime} the computer running times with conforming and nonconforming time grids on a fixed mesh $h=1/32$. We see that the use of nonconforming time grids is efficient in terms of accuarcy and computational cost.

\begin{table}[h!]
\begin{center}
\begin{tabular}{| c ||  c | c | c | c | }
\firsthline
 \multirow{2}{*}{$\Delta t$}  &  \multicolumn{2}{c|}{Linear viscosities}  & \multicolumn{2}{c|}{Nonlinear viscosities}   \\[1pt] \cline{2-5}
    & Conforming & Nonconforming & Conforming & Nonconforming \\[1pt]
\hline
$0.8$  		&  72 	&   & 167  &      \\[1pt]   & & 83	&  & 188  \\[1pt]  
$0.4 $ 		&  154  &   &  348 &    \\[1pt]     & & 180	& 	&	420  \\[1pt]  
$0.2 $  		&  310 	&  &  727  &    \\[1pt]    & & 351	& 	&	791 \\[1pt]  
$0.1 $ 		&  632	&  &  1447  &   \\[1pt]     & & 697	& 	&	1646  \\[1pt]  
$0.05 $ 	& 1262  &  &  2791 &   \\
\lasthline
\end{tabular}
\end{center}
\caption{Comparison of the computer running times  (in seconds) of conforming and nonconforming time grids, with MINI elements for the Stokes problem and Raviart-Thomas elements for the Darcy problem on a fixed mesh $h=1/32$. Note that for the nonlinear viscosities, two Gauss-Newton iterations are performed.} \label{tab:RTcputime}
\end{table}

%\begin{table}[h!]
%\begin{center}
%\begin{tabular}{| c ||  c | c | c | c | }
%%
%\firsthline
%%
% \multirow{2}{*}{$\Delta t$}  &  \multicolumn{2}{c}{Linear viscosities}  & \multicolumn{2}{|c|}{Nonlinear viscosities}   \\[1pt] \cline{2-5}
%    & Conforming & Nonconforming & Conforming & Nonconforming \\[1pt]
%\hline
%$0.8$  		&  1150 	    &   &  3402&      \\[1pt]   & & 1332 &  & 3745  \\[1pt]  
%$0.4 $ 		&  2424    &   & 5561 &    \\[1pt]     & &  2673 & 	&	6716  \\[1pt]  
%$0.2 $  		&  4791 	&  &  11058  &    \\[1pt]    & & 5506	& 	&	12261  \\[1pt]  
%$0.1 $ 		&  9763 	&  &  21690 &   \\[1pt]     & & 10449 & 	&	23963  \\[1pt]  
%$0.05 $ 	&  19636  &  & 42101 &   \\
%\lasthline
%\end{tabular}
%\end{center}
%\caption{Comparison of computation time  (in seconds) of conforming and nonconforming time grids, with MINI elements for the Stokes problem and Raviart-Thomas elements for the Darcy problem on a fixed mesh $h=1/64$. Note that for the nonlinear viscosities, two Gauss-Newton iterations are performed.}
%\end{table}
%
%
%  TEST CASE 2 (Yotov's paper)
%
%
\subsection{Test case 2: flow driven by a pressure drop} \label{subsec:t2}
\begin{figure}[H]
    \centering
    \begin{subfigure}[b]{0.45\textwidth}
        \includegraphics[scale=0.35]{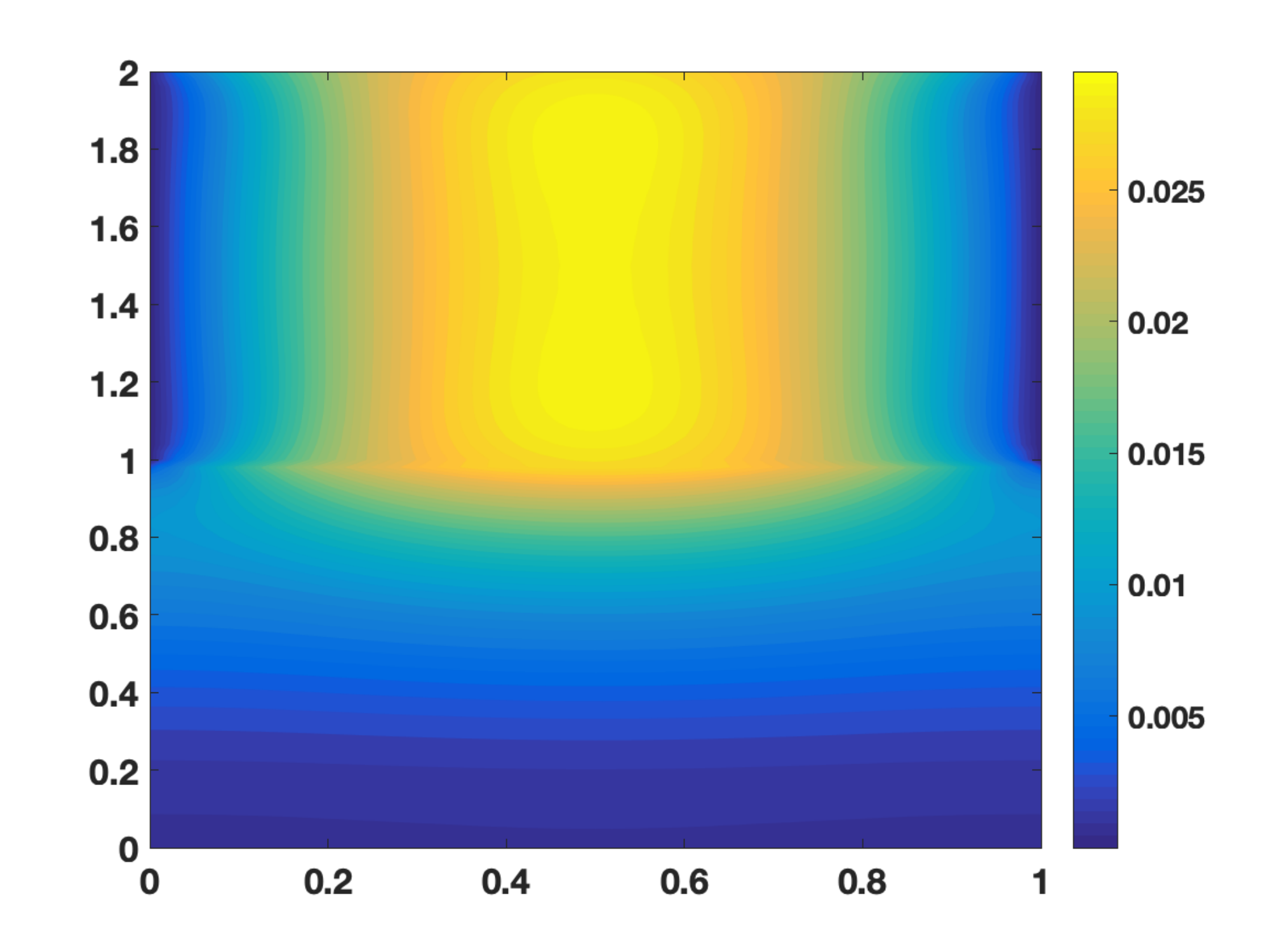}
    \end{subfigure}
    \begin{subfigure}[b]{0.45\textwidth}
        \includegraphics[scale=0.35]{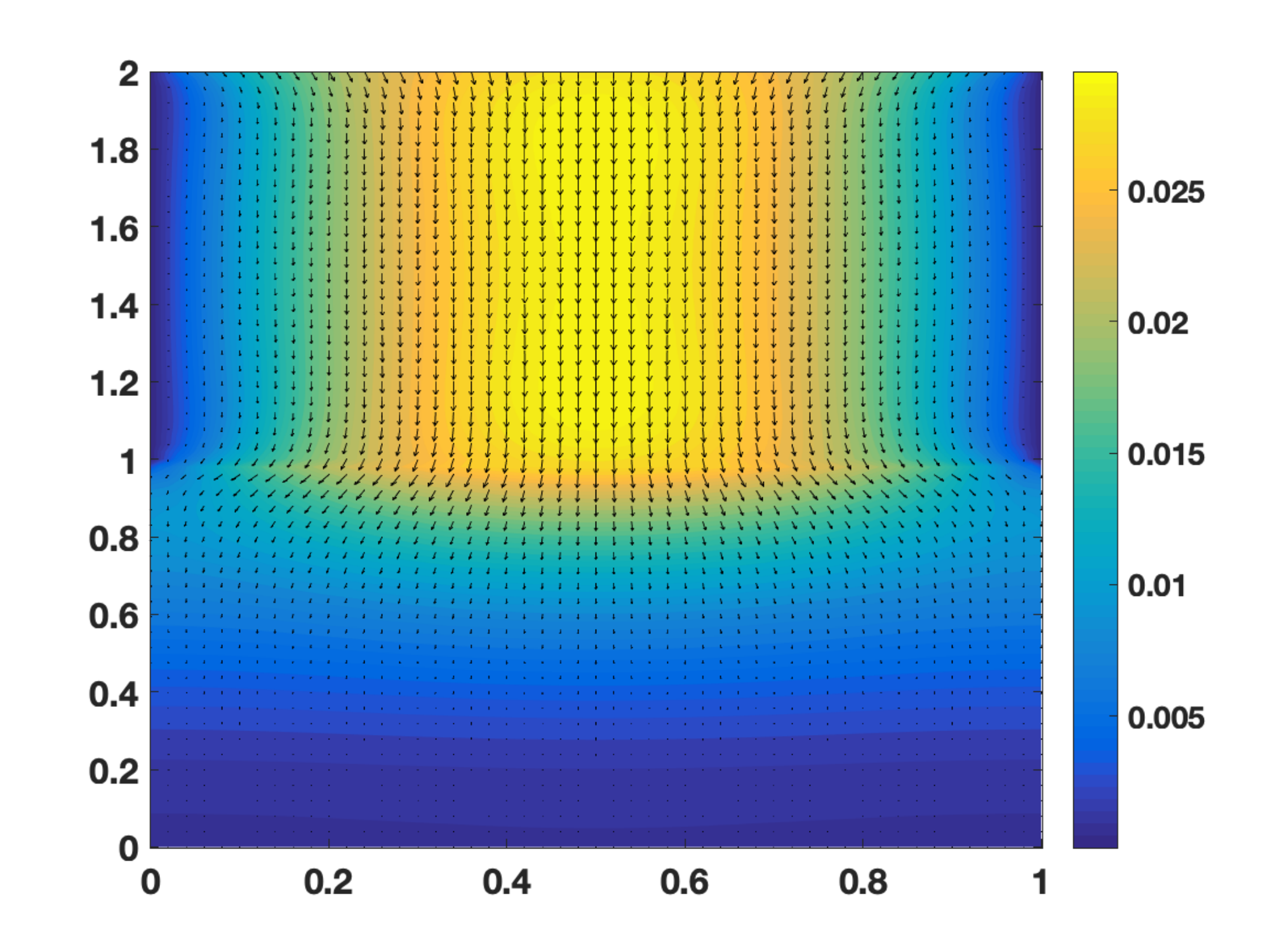}
    \end{subfigure} 
    \caption{[Test case 2] Velocity magnitude and velocity vector at $T=1$.}\label{fig:velo}
\end{figure}

In this test case, the flow is driven by a pressure drop: on the top boundary of $\Omega_{f}$ we set $p_{\text{in}}=1$ and on the bottom boundary of $\Omega_{p}$, $p_{\text{out}}=0$, which is also chosen as the initial condition for the Darcy pressure. Along the left and right boundaries, we impose no-slip boundary condition for the Stokes flow and no-flow boundary condition for the Darcy flow. We also set zero velocity initial condition for the Stokes problem. The parameters are $\kappa=1$, $K_{f}=K_{p}~=~1$, $\nu_{f, \infty}=\nu_{p, \infty}=1$, $\nu_{f, 0}=\nu_{p,0}=10$, $r_{f}=r_{p}=1.35$ and $\alpha=1$. The simulation time is $T=1$. For this test case, we use cell conservative spatial discretization, i.e. MINI elements for the Stokes flow and the Raviart-Thomas of order 1 elements for the Darcy flow. The velocity magnitude and vector at the final time are shown Figure~\ref{fig:velo}. 

We compute the reference solution on a mesh size $h=1/32$ and $\Delta t_{\text{ref}}=0.01$. We want to verify the convergence in time of the global-in-time domain decomposition method with nonconforming time grids: $\Delta t_{f}=\sfrac{\Delta t_{p}}{2}$. Table~\ref{tab:testcase2} shows the errors of the nonlinear Stokes and Darcy problems at $T=1$ with a fixed mesh size $h=1/32$, first order convergence in time is observed. In Tables~\ref{tab:t2Stokes} and \ref{tab:t2Darcy}, we compare the accuracy in time of the conforming and nonconforming time grids. In particular, the errors (with nonconforming time grids) in the fluid domain are close to those with fine conforming time steps, while those in the porous medium are close to those with coarse conforming time steps.  

\begin{table}[http!]
	\setlength{\extrarowheight}{4pt}
	\centering
	\small
\begin{tabular}{| l | l | l  l | l  l | l  l | l  l |} \hline 
		\multicolumn{2}{|c|}{Time step} & \multicolumn{2}{c|}{$\bu_{f}$} & \multicolumn{2}{c|}{$p_{f}$} & \multicolumn{2}{c|}{$\bu_{p}$} & \multicolumn{2}{c|}{$p_{p}$}\\ \hline
		$\Delta t_{f}$& $\Delta t_{p}$ & \multicolumn{2}{c|}{$H^{1}$ error}  & \multicolumn{2}{c|}{$L^{2}$ error} & \multicolumn{2}{c|}{$H^{\text{div}}$ error}& \multicolumn{2}{c|}{$L^{2}$ error} \\ \hline
		$1/4$ & $1/2$ &   3.44e-04 &        & 2.80e-03 &         &   4.56e-03 &      & 2.17e-03 &\\  
		$1/8$ & $1/4$ &  1.49e-04  & [1.21] & 1.37e-03  & [1.03] & 2.18e-03  & [1.07] & 1.12e-03  & [0.95]\\ 
		$1/16$ & $1/8$ &  5.60e-05 & [1.41] & 6.51e-04 & [1.07]  &  1.04e-03 & [1.07] & 5.48e-04 & [1.03]\\ 
		$1/32$ & $1/16$ & 1.63e-05 & [1.78] & 2.95e-04 & [1.14] & 4.70e-04 & [1.15]  & 2.53e-04 & [1.12]\\  \hline
	\end{tabular} 
\caption{[Test case 2] Errors for the nonlinear Stokes and Darcy problems at $T=1$ with a fixed mesh size $h=1/32$. }  \label{tab:testcase2}
\end{table}

\begin{table}[http!]
	\small
	\setlength{\extrarowheight}{4pt}
	\centering
	\begin{tabular}{| l| c |c| c | c | c |  } \hline  
		\multirow{2}{*}{Time grids} & \multirow{2}{*}{$\Delta t_{f}$}& \multirow{2}{*}{$\Delta t_{p}$} & \multicolumn{2}{c|}{$\bu_{f}$} & $p_{f}$ \\ \cline{4-6}
		 & &  & $L^{2}$ error  & $H^{1}$ error & $L^{2}$ error \\ \hline
		Conforming coarse &$1/8$ & $1/8$ & 3.01e-05  & 1.03e-04  & 6.71e-04 \\ \hline 
		Nonconforming & $1/16$ & $1/8$ &  1.64e-05 & 5.60e-05 & 6.51e-04  \\ \hline  
		Conforming fine & $1/16$ & $1/16$ & 1.38e-05 & 4.68e-05 & 3.08e-04  \\ \hline 
	\end{tabular}  
\caption{[Test case 2] Errors for the nonlinear Stokes problem at $T=1$ with a fixed mesh size $h=1/32$. }  \label{tab:t2Stokes}
\end{table}	
	
\begin{table}[http!]
	\small
	\setlength{\extrarowheight}{4pt}
	\centering
		\begin{tabular}{|l | c |c| c | c | c |  } \hline 
		\multirow{2}{*}{Time grids} & \multirow{2}{*}{$\Delta t_{f}$}& \multirow{2}{*}{$\Delta t_{p}$} & \multicolumn{2}{c|}{$\bu_{f}$} & $p_{f}$ \\ \cline{4-6}
		& &  & $L^{2}$ error  & $H^{\text{div}}$ error & $L^{2}$ error \\ \hline
		Conforming coarse & $1/8$ & $1/8$ &  2.11e-04 &  1.05e-03 & 5.56e-04 \\ \hline 
		Nonconforming & $1/16$ & $1/8$ & 2.08e-04 & 1.04e-03  &  5.48e-04  \\ \hline 
		Conforming fine & $1/16$ & $1/16$ & 9.71-05 & 4.78e-04  & 2.59e-04  \\ \hline 
	\end{tabular}  
\caption{[Test case 2] Errors for the nonlinear Darcy problem at $T=1$ with a fixed mesh size $h=1/32$. }  \label{tab:t2Darcy}
\end{table}	

\begin{figure}[H]
    \centering
\includegraphics[scale=0.4]{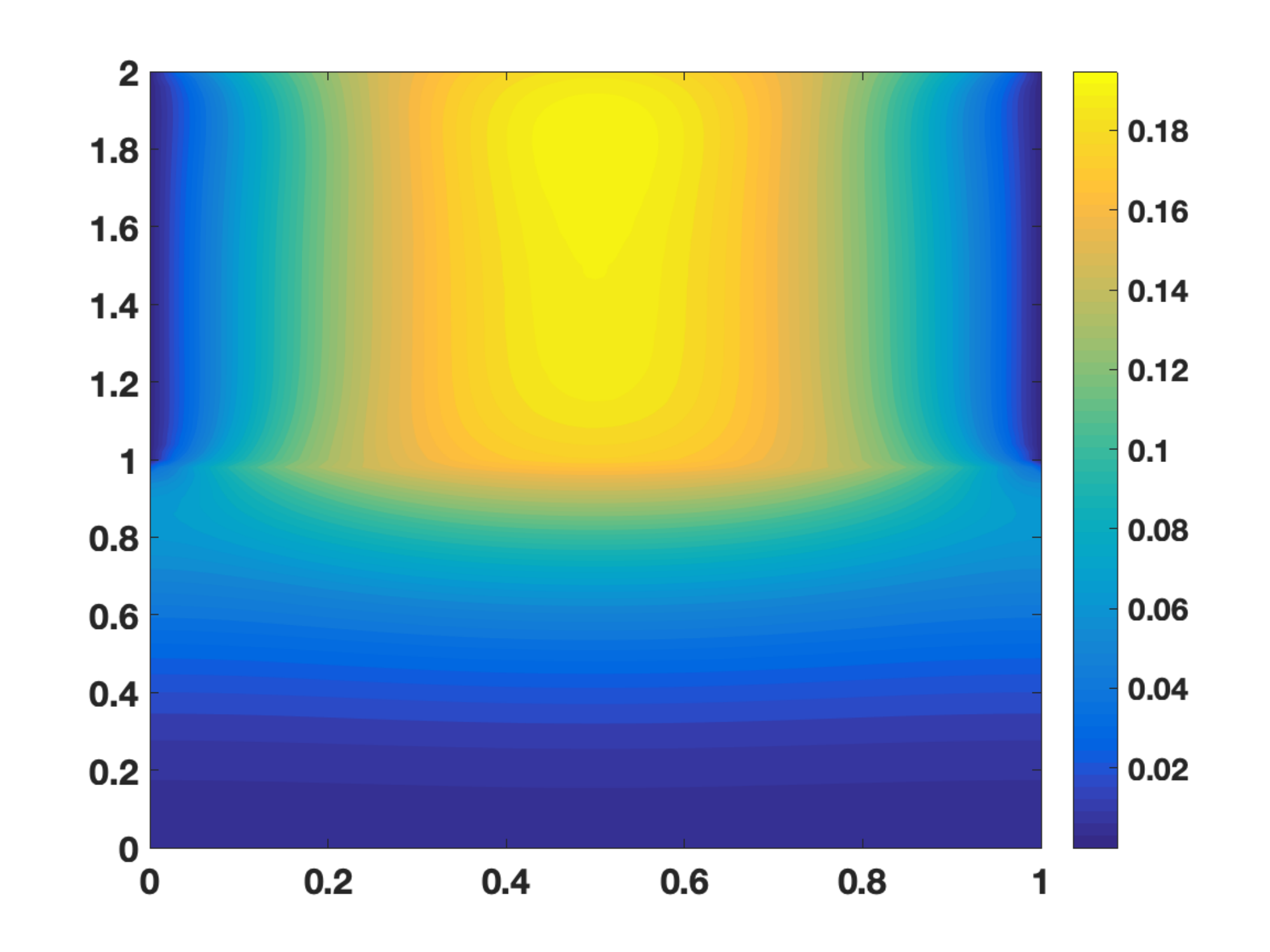}
\caption{[Test case 2 with discontinuous paramters] Velocity magnitude at $T=1$.}\label{fig:velodis}
\end{figure}

Next, we consider the case with discontinuous parameters. In particular, we have, for the Stokes problem, $K_{f}=1, \, \nu_{f, \infty}=0.5, \, \nu_{f, 0}=1,$ and for the Darcy problem, $K_{p}~=0.001, \, \nu_{p, \infty}=1, \, \nu_{p,0}=10.$ As before, we impose smaller time step in the fluid region and larger time step in the porous medium: $\Delta t_{f} =  \Delta t_{p}/2=0.125$.
The velocity magnitude at $T=1$ is depicted in Figure~\ref{fig:velodis} and the errors are reported in Table~\ref{tab:t2dis}, which again shows that the accuracy in time is well-preserved with nonconforming time grids.

\begin{table}[http!]
	\small
	\setlength{\extrarowheight}{4pt}
	\centering
\begin{tabular}{|l | c| c| c |c| c | c |  } \hline 
		Time grids & $\Delta t_{f}$ & $\Delta t_{p}$ & $u_{f}$ & $p_{f}$ & $u_{p}$ & $p_{p}$ \\ \hline
		Conforming coarse &  $1/4$ & $1/4$ & 1.48e-03 & 2.02e-02 & 5.41e-03 & 1.45e-02 \\ \hline 
		Nonconforming  &  $1/8$ & $1/4$  & 2.33e-04 & 1.69e-02  &  4.70e-03  & 1.32e-02 \\ \hline 
		Conforming fine  &  $1/8$ & $1/8$ &  2.27e-04 & 1.00e-02 & 2.77e-03 & 7.63e-03 \\ \hline 
	\end{tabular} \vspace{4pt}\\
\caption{[Test case 2 with discontinuous parameters] $L^{2}-$errors for the nonlinear Stokes and Darcy problems at $T=1$ with $h=1/32$. }  \label{tab:t2dis}
\end{table}

\section{Conclusion} \label{sec:conclu}
We have introduced a decoupling scheme for the nonlinear Stokes-Darcy system, based on the time-dependent interface operators. The scheme is an implicit type that requires iterations between subdomains; the subproblems, time-dependent Stokes and Darcy equations, are solved using local time-stepping algorithms, respectively. The space-time domain decomposition method allows us to independently solve each subproblem using existing local solvers and enables the use of nonconforming time grids as well as different time-stepping algorithms for local problems. 
For numerical tests of the proposed algorithm two numerical examples were considered; the first is a non-physical problem with the known exact solution and the second is a flow problem driven by a pressure drop. Numerical results confirm that the algorithm simulates the model problem at the optimal order of accuracy and its efficiency is improved with the use of nonconforming time grids and the preconditioner for GMRES iterations. Although the model system is nonlinear, only one or two Newton iterations were needed within the given tolerance range, yielding the optimal accuracy in our test cases. 

Some future directions for this work include extending the approach to more complex coupled problems such as the coupled Stokes-Darcy system with transport and a fluid flow coupled with a quasi-static poroelastic medium. In particular, because of the use of local time stepping, we expect that this approach is efficiently applicable to multiphysics problems, where local problems are in different time scales, e.g., 
fluid flows interacting with clays or soils. Many such examples are found in applications of geomechanics and the quasi-static Biot's consolidation model \cite{Biot} is often considered for a deformable porous medium. 
In the Biot model, the fluid motion in the porous medium is described by Darcy's law, while the deformation of the medium is governed by the linear elasticity. Interface conditions for the (Navier-)Stokes-Biot system are more complex than those of the Stokes-Darcy system, however, we expect that a similar approach can be considered for the large multiphysics problem to be turned into a time dependent Steklov-Poincar\'e operator equation.  
We are also currently investigating an optimized Schwartz waveform relaxation (OSWR) method using Robin transmission conditions for the Stokes-Darcy model considered in this work. Details concerning the development, analysis and numerical implementation of space-time domain decomposition based on OSWR is a subject of a forthcoming paper.

\bibliographystyle{spmpsci}

\end{document}